\documentclass[11pt,english,a4paper]{smfart}
\usepackage[english]{babel}
\usepackage{xcolor}
\usepackage{pbsi}
\usepackage[T1]{fontenc}
\usepackage{stmaryrd}
\usepackage{mathabx}
\usepackage[mathscr]{euscript}
\usepackage[pagebackref,colorlinks=true,linkcolor=blue,citecolor=blue,urlcolor=red, hypertexnames=true]{hyperref}

 \usepackage[abbrev,backrefs,nobysame]{amsrefs}
 \usepackage{caption}
 \usepackage{enumerate}

\usepackage{amsmath,amssymb,bm,amsfonts}
\usepackage{nicematrix}
\usepackage[all]{xy}
\entrymodifiers={+!!<0pt,\fontdimen22\textfont2>}

\usepackage{mathrsfs}

\usepackage{graphicx}

\usepackage{color}

\newtheorem{theorem}{Theorem}[section]
\newtheorem{lemma}[theorem]{Lemma}

\newtheorem{corollary}[theorem]{Corollary}
\newtheorem{proposition}[theorem]{Proposition}

\newtheorem*{remark*}{\it Remark}

\newcommand{\flba}[2]{
\xymatrix@C15pt{#1\ar@{|->}[r]&#2}}
\newcommand{\flcourte}[2]{
\xymatrix@C12pt{#1\ar[r]&#2}}

{\theoremstyle{definition}}

{\theoremstyle{definition}}

\theoremstyle{definition}

\newtheorem{question}[theorem]{Question}
\newtheorem{fact}[theorem]{Fact}

{\theoremstyle{definition}\newtheorem{remark}[theorem]{Remark}}

\def\D{\ensuremath{\mathbb D}}

\def\T{\ensuremath{\mathbb T}}
\def\R{\ensuremath{\mathbb R}}
\def\Z{\ensuremath{\mathbb Z}}

\def\C{\ensuremath{\mathbb C}}
\def\Q{\ensuremath{\mathbb Q}}
\def\N{\ensuremath{\mathbb N}}

\def\u{{{\mathcal U}}}

\def\en{E_{N}}
\def\phip{{{\Phi_+(X)}}}

\def\bth{\begin{theorem}}
\def\blm{\begin{lemma}}
\def\bpr{\begin{proposition}}
\def\bpf{\begin{proof}}
\def\epf{\end{proof}}
\def\epr{\end{proposition}}
\def\elm{\end{lemma}}
\def\eth{\end{theorem}}
\def\bco{\begin{corollary}}
\def\eco{\end{corollary}}
\def\be{\begin{enumerate}}
\def\ee{\end{enumerate}}
\def\bea{\begin{enumerate}[\rm (a)]}
\def\beun{\begin{enumerate}[\rm (1)]}
\def\bei{\begin{enumerate}[\rm (i)]}

\newcommand{\pss}[2]{\ensuremath{{\langle #1,#2\rangle}}}

\newcommand{\vertiii}[1]{{\left\vert\kern-0.25ex\left\vert\kern-0.25ex\left\vert #1 
    \right\vert\kern-0.25ex\right\vert\kern-0.25ex\right\vert}}

\newcommand{\bmx}{{\mathcal{B}}_{M}(X)}
\newcommand{\pmx}{{\mathcal{P}}_{M}(X)}

\newcommand{\bbx}{{\mathcal{B}}_{1}(X)}
\newcommand{\ppx}{{\mathcal{P}}_{1}(X)}

\newcommand{\wot}{\texttt{WOT}}

\newcommand{\sot}{\texttt{SOT}}
\newcommand{\sote}{\texttt{SOT}\mbox{$^{*}$}}
\newcommand{\bx}{{\mathcal B}(X)}
\newcommand{\px}{{\mathcal P}(X)}

\newcommand{\lan}{\langle}
\newcommand{\ran}{\rangle}
\newcommand{\akl}{a_{k,l}}

\numberwithin{equation}{section}

\author[V. Gillet]{Valentin Gillet}
\address[V. Gillet]{Laboratoire Paul Painlev\'e, UMR 8524\\
Universit\'{e} de Lille\\
Cit\'e Scientifique, B\^atiment M2\\
59655 Villeneuve d'Ascq Cedex 
(France)}
\email{valentin.gillet@univ-lille.fr}

\thanks{This work was supported in part by the project COMOP of the French National Research Agency (grant ANR-24-CE40-0892-01) and by the Labex CEMPI (ANR-11-LABX-0007-01). The author also acknowledges the support of the CDP C2EMPI, as well as the French State under the France-2030 programme, the University of Lille, the Initiative of Excellence of the University of Lille, and the European Metropolis of Lille for their funding and support of the R-CDP-24-004-C2EMPI project.}

\begin{document}

\title[Typical $\ell_p\,$-$\,$space contractions]{Typical properties of positive contractions and the invariant subspace problem}


\begin{abstract} 
In this paper, we first study some elementary properties of a typical positive contraction on $\ell_q$ for the \sot\, and the \sote \,topologies. Using these properties, we prove that a typical positive contraction on $\ell_1$ (resp. on $\ell_2$) has a non-trivial invariant subspace for the \sot\, topology (resp. the \sot\, and the \sote \,topologies). We then focus on the case where $X$ is a Banach space with a basis. We prove that a typical positive contraction on a Banach space with an unconditional basis has no non-trivial closed invariant ideals for the \sot\, and the \sote \,topologies. In particular, this shows that when $X = \ell_q$ with $1 \leq q < \infty$, a typical positive contraction $T \in (\ppx,\sot)$ (resp. $T \in (\ppx, \sote)$ when $1 < q < \infty$) does not satisfy the Abramovich, Aliprantis and Burkinshaw criterion, that is, there is no non-zero positive operator in the commutant of $T$ which is quasinilpotent at a non-zero positive vector of $X$. Finally, we prove that, for the \sote\, topology, a typical positive contraction on a reflexive Banach space with a monotone basis does not satisfy the Abramovich, Aliprantis and Burkinshaw criterion.
\end{abstract}

\maketitle

\par\bigskip
\section{Introduction}\label{Section introduction}

Throughout this paper, $(X, \lVert . \rVert)$ will be a complex separable infinite-dimensional Banach space and $H$ will be a complex separable infinite-dimensional Hilbert space. The space of all bounded operators on $X$ equipped with the operator norm $\lVert . \rVert$ will be denoted by $\bx$.

If $M>0$, we denote by $\bmx$ the closed ball of radius M of $\bx$. With this notation the set $\bbx$ is just the set of contractions on $X$. If $(e_n)_{n \geq 0}$ is a basis of $X$, we say that a vector $x = \displaystyle \sum_{n\geq 0} x_n e_n$ of $X$ is positive whenever $x_n \geq 0$ for all $n \geq 0$. We write $x \geq 0$ when $x$ is a positive vector of $X$. We say that an operator $T \in \bx$ is positive whenever $Tx \geq 0$ for every $x \geq 0$. The space of all bounded positive operators on $X$ will be denoted by $\px$, and we denote by $\pmx := \bmx \cap \px$ the positive operators in $\bmx$.

\par\smallskip

Given a topology $\tau$ on $\mathcal{Y} = \bmx$ or $\mathcal{Y} = \pmx$ that turns $\mathcal{Y}$ into a Baire space, we say that a property (P) of operators on $X$ is typical for the topology $\tau$ if the set $\{ T \in \mathcal{Y} : T \; \textrm{satisfies the property (P)}\}$ is comeager in $\mathcal{Y}$, that is, contains a dense $G_\delta$ subset of $\mathcal{Y}$. A property (P) of operators on $X$ is atypical for the topology $\tau$ if its negation is typical for $\tau$. Since the space $(\bbx, \lVert . \rVert)$ is usually not separable and so is not Polish, we will be focusing on the Strong Operator Topology and on the Strong* Operator Topology. Recall that the Strong Operator Topology is the topology on $\bx$ defined by the seminorms $\lVert . \rVert_x$, where $\lVert T \rVert_x = \lVert Tx \rVert$ for every $x \in X$. The Strong* Operator Topology is the topology on $\bx$ defined by the seminorms $\lVert . \rVert_x$ and $\lVert . \rVert_{x^*}$, where $\lVert T \rVert_x = \lVert Tx \rVert$ and $\lVert T \rVert_{x^*} = \lVert T^* x^* \rVert$ for every $x \in X$ and every $x^* \in X^*$. If $(T_i)$ is a net in $\bx$ and $T \in \bx$, we have
$$
T_i \underset{i}{\overset{\sot}{\longrightarrow}}T \iff T_i x \underset{i}{\overset{\lVert . \rVert}{\longrightarrow}}Tx \quad \textrm{for every} \; x \in X,
$$
and
$$
T_i \underset{i}{\overset{\sote}{\longrightarrow}}T \iff
\left\{
\begin{array}{ll}
         T_i x \underset{i}{\overset{\lVert . \rVert}{\longrightarrow}}Tx   &\textrm{for every} \; x \in X,\\
        T^*_i x^* \underset{i}{\overset{\lVert . \rVert}{\longrightarrow}}T^* x^*  &\textrm{for every} \; x^* \in X^*.
    \end{array}
\right.
$$
We now denote by \sot \, (resp. by \sote) the Strong Operator Topology (resp. the Strong* Operator Topology) on $\bx$.
For every $M > 0$, the closed ball $(\bmx, \sot)$ is Polish and when $X^*$ is separable, the space $(\bmx, \sote)$ is Polish.

\par\smallskip The notion of a typical property of contractions was initiated by Eisner in \cite{E}. She proved in \cite{E} that a typical contraction on a separable Hilbert space is unitary for the Weak Operator Topology. This notion was studied in more depth by Eisner and Mátrai for operators on a separable Hilbert space for other topologies as, for example, the \sot \, and the \sote \, topologies. It is proved in \cite{EM} that if $H$ is a separable Hilbert space, a typical $T \in (\mathcal{B}_1(H), \sot)$ is unitarily equivalent to the infinite-dimensional backward unilateral shift operator on $\ell_2(\Z_+\times \Z_+)$ and in particular has a non-trivial invariant subspace, that is, a closed subspace $Y \subseteq X$ such that $Y \ne \{0 \}, Y \ne X$ and $T(Y) \subseteq Y$. Grivaux, Matheron and Menet have also studied typical properties of contractions on $\ell_q$-spaces (\cite{GMM}, \cite{GMM1} and \cite{GMM2}). Their initial goal was to determine whether a typical $T \in (\mathcal{B}_1(\ell_q), \sot)$ or $T \in (\mathcal{B}_1(\ell_q), \sote)$ has a non-trivial invariant subspace.

\par\smallskip
This motivation comes from the famous invariant subspace problem, which can be stated as follows: given an infinite-dimensional separable complex Banach space $Z$, does every bounded operator $T \in \mathcal{B}(Z)$ have a non-trivial invariant subspace? The answer to this problem is negative for some non-reflexive Banach spaces: Enflo constructed a non-reflexive Banach space which supports an operator without non-trivial invariant subspaces (see \cite{En}), and Read constructed operators on $\ell_1$ and on $c_0$ (called Read's operators) without non-trivial invariant subspaces (see \cite{R2} and \cite{R3}). The invariant subspace problem still remains open for reflexive Banach spaces and in particular for separable Hilbert spaces. However, the problem has a positive answer for some classes of operators: for example, Lomonosov proved in \cite{L} that if an operator $T$ contains a non-scalar operator in its commutant that commutes with a non-zero compact operator, then it has a non-trivial invariant subspace; Brown, Chevreau and Pearcy proved in \cite{BCP2} that every contraction on a Hilbert space whose spectrum contains the unit circle has a non-trivial invariant subspace. The invariant subspace problem can be restricted to positive operators as follows: given an infinite-dimensional complex Banach space $Z$ with a basis, does every bounded positive operator on $Z$ have a non-trivial invariant subspace? The answer to this problem is still unknown. Moreover, even if Read's operators give a counter-example to the invariant subspace problem on $Z = \ell_1$, it was proved in \cite{Troi} that the modulus of a Read's operator on $\ell_1$ has a positive eigenvector and so the invariant subspace problem for positive operators still remains open for the case $Z = \ell_1$. A major result regarding the existence of a non-trivial invariant subspace for positive operators, due to Abramovich, Aliprantis and Burkinshaw (\cite[Theorem 2.2]{AAB}), is the following.

\smallskip

\bth \label{thintropos}
    Let $X$ be a Banach space with a basis and $T$ be a positive operator on $X$. If there exists a non-zero positive operator $A$ on $X$ which is quasinilpotent at a certain non-zero positive vector of $X$ and such that $A T = T A$, then $T$ has a non-trivial invariant subspace. 
\eth

\smallskip

In particular, Theorem \ref{thintropos} implies the following corollary.
\smallskip

\bco
Let $X$ be a Banach space with a basis. Every positive operator on $X$ which is quasinilpotent at a non-zero positive vector of $X$ has a non-trivial invariant subspace.
\eco

\smallskip

For more details on the invariant subspace problem and on known results on this subject, we refer to  \cite{CE}, \cite{CP} and \cite{RR}.

\par\smallskip
The aim of this article is to study the invariant subspace problem from the point of view of typicality, that is, we are interested in the question of whether a typical $T \in (\ppx, \tau)$ has a non-trivial invariant subspace, when $\tau = \sot$ or $\tau = \sote$. We will mostly focus on the case $X = \ell_q$ with $1 \leq q < \infty$, but we will also generalize some results to Banach spaces with an unconditional basis or with a monotone basis. In the case where $X = \ell_q$, the basis $(e_n)_{n \geq 0}$ of $X$ will be the canonical basis, where we recall that $e_n(k) = 1$ if $n = k$ and $e_n(k) = 0$ if $n \ne k$, for every $n,k \geq 0$.

\subsection{Notations}\label{notation}

We introduce here some notation that will be used throughout the article.
 
\par\smallskip
- We denote by $\Q_{>0}$ the set of positive rational numbers.

\smallskip

- The open unit disk of $\C$ will be denoted by $\D$ and the unit circle of $\C$ will be denoted by $\T$.

\smallskip
- If $Z$ is a Banach space, we denote by $S_Z$ the unit sphere of $Z$.

\smallskip
- The closed linear span of $(x_i)_{i \in I} \subseteq X$ will be written as $[x_i : i \in I]$. 

\smallskip
- If $X$ is a Banach space with a basis $(e_n)_{n \geq 0}$, we denote by $E_N$ the subspace $[e_0,...,e_N ]$ for every $N \geq 0$ and by $F_N$ the subspace $ [e_j : j > N]$ for every $N \geq 0$.

\smallskip
- If $X$ is a Banach space with a basis $(e_n)_{n \geq 0}$, we denote by $P_N$ the canonical projection onto $E_N$ and by $Q_N$ the canonical projection onto $F_N$. The biorthogonal functionals sequence associated to the basis $(e_n)_{n \geq 0}$ will be denoted by $(e_n^*)_{n \geq 0}$.

\smallskip
- If $T \in \bx$, we respectively write $\sigma(T), \, \sigma_{ap}(T), \, \sigma_p(T)$ and $\sigma_{ess}(T)$ for the spectrum, the approximate spectrum, the point spectrum and the essential spectrum of $T$.

\smallskip

\subsection{Main results} \label{mainresults}

    We start by presenting in Section \ref{sectiontool} some tools that will be useful throughout the article. We properly define the notion of a positive operator on a Banach space with a basis. In particular, Proposition \ref{polishsotsote} states that the set of positive contractions on a Banach space with a basis is Polish for both the \sot\, and the \sote \,topologies, providing the necessary setting for the study of typical properties of positive contractions.

\smallskip
    
   Theorem \ref{sotsimsoteintro} connects the comeager sets of $\ppx$ for the topologies \sot\, and \sote\, in the case where $X = \ell_q$ with $q > 2$.

    \smallskip
    
   In Section \ref{Section3}, we put together some elementary properties of a typical positive contraction for the \sot \, and the \sote \, topologies. Proposition \ref{normspec} will play an important role for Section \ref{Section4} in order to study the invariant subspace problem for a typical positive contraction on $\ell_2$ and on $\ell_1$.

\smallskip

We also study the point spectrum of a typical positive contraction on $\ell_q$ (Corollaries \ref{corvalprop} and \ref{corvalpropbis}). 

\smallskip

Eisner and Mátrai proved in \cite{EM} that a typical $T \in (\mathcal{B}_1(\ell_2), \sot)$ is unitarily equivalent to the infinite-dimensional backward unilateral shift operator on $\ell_2(\Z_+ \times \Z_+)$ and in particular, this implies that a typical contraction on $\ell_2$ has a non-trivial invariant subspace. To do so, they proved first that a typical contraction $T \in (\mathcal{B}_1(\ell_2), \sot)$ is such that $T^*$ is an isometry. It turns out that this is no longer the case for positive contractions, that is, the adjoint of a typical positive contraction on $\ell_2$ is not an isometry for the \sot \, topology. 

\smallskip

\bpr \label{introisometry}
Let $X = \ell_q$ with $1 < q < \infty$. A typical $T \in (\ppx, \emph{\sot})$ (resp. $T \in (\ppx, \emph{\sote})$) is such that $T^*$ is not an isometry. 
\epr

\smallskip
The case $X = \ell_1$ is very different from the case $X = \ell_q$ with $q > 1$, as the following proposition shows.

\smallskip

\bpr
If $X = \ell_1$, then a typical $T \in (\mathcal{P}_1(X),\emph{\sot})$ is such that $T^* $ is a non-surjective isometry and such that $T-\lambda$ is surjective for every $\lambda \in \D$.
\epr

\smallskip

Proposition \ref{introisometry} shows in particular that a property can be typical in $\ppx$ but can be atypical in $\bbx$ for a certain topology $\tau$ on $\bbx$.

\smallskip

In Section \ref{Section4}, we first explain why a typical positive contraction on $\ell_1$ and on $\ell_2$ has a non-trivial invariant subspace, and then we focus on the case where $X$ is a Banach space with a basis. In particular, it includes the case $X = \ell_q$ with $1 < q \ne 2 < \infty$. 

We say that a positive operator on a Banach space $X$ with a basis satisfies the Abramovich, Aliprantis and Burkinshaw criterion (abbreviated AAB criterion) if this operator satisfies the hypotheses of Theorem \ref{thintropos}. Any positive operator on a Banach space $X$ with a basis satisfying the AAB criterion has a non-trivial invariant subspace and when $X = \ell_q$ with $1 \leq q <\infty$, any operator satisfying the AAB criterion has a non-trivial closed invariant ideal (\cite[Theorem 2.2]{AAB2}), that is, a closed vector subspace $V$ of $X$ such that $\lvert x \rvert \leq \lvert y \rvert$ and $y \in V$ imply $x \in V$, for every $x,y \in X$. Here, $\lvert x \rvert$ is the positive vector of $X$ whose coordinates are the modulus of the coordinates of the vector $x \in X$. Troitsky and Radjavi gave the following characterization of positive operators on a Banach space with an unconditional basis admitting a non-trivial invariant closed ideal (\cite[Proposition 1.2]{RT}).

\smallskip

\bpr \label{critTroitsky}
Let $X$ be a Banach space with an unconditional basis $(e_n)_{n \geq 0}$ and let $T$ be a positive operator on $X$.
The operator $T$ has no non-trivial closed invariant ideals if and only if the following property holds:
$$
\forall \, i \ne j \in \Z_+, \, \exists n \in \Z_+ : \, \lan e_j^*, T^n e_i \ran > 0.
$$
\epr

\smallskip

Thanks to this characterization, we obtain the following result.

\smallskip

\bpr \label{propidealinv}
Let $X$ be a Banach space with an unconditional basis. A typical $T \in (\ppx, \emph{\sot})$ (resp. $T \in (\ppx, \emph{\sote})$ when $X^*$ is separable) has no non-trivial closed invariant ideals.  
\epr

\smallskip

In particular, Proposition \ref{propidealinv} implies the following corollary.

\smallskip

\bco \label{corAABincond}
Let $X = \ell_q$ with $1 \leq q < \infty$. A typical $T \in (\ppx, \emph{\sot})$ (resp. $T \in (\ppx, \emph{\sote})$ when $1 < q < \infty$) does not satisfy the AAB criterion. 
\eco

\smallskip

The main result of this article is the following generalization of Corollary \ref{corAABincond} for the \sote \,topology to Banach spaces with a monotone basis. 

\smallskip

\bth \label{thAABmonotone}
Let $X$ be a reflexive Banach space with a monotone basis. A typical $T \in (\ppx, \emph{\sote})$ does not satisfy the AAB criterion.
\eth

\smallskip

Finally, we end Section \ref{Section4} with the following result.

\smallskip

\bco \label{corqny}
Let $X$ be a Banach space with a basis. A typical $T \in (\ppx, \emph{\sot})$ (resp. $T \in (\ppx, \emph{\sote})$ when $X^*$ is separable) is not quasinilpotent at any non-zero positive vector of $X$. 
\eco

\smallskip

These results highlight the fact that the existence of invariant subspaces for typical operators is a delicate matter, even in the restricted setting of positive contractions.

\section{Useful tools} \label{sectiontool}

\subsection{Positive operators on $X$.} \label{posop}

    If $X$ is a Banach space with a basis $(e_n)_{n \geq 0}$, we denote by $\mathcal{C}^+$ the positive cone of $X$, that is:
    \begin{align*}
        \mathcal{C}^{+} := \{ x = \displaystyle \sum_{n \geq 0} x_n e_n : x_n \geq 0, \; \textrm{for every} \; n \geq 0  \}.
    \end{align*}

    We have $\mathcal{C}^+ + \mathcal{C}^+ \subseteq \mathcal{C}^+$, $\alpha \mathcal{C}^+ \subseteq \mathcal{C}^+$ for every $\alpha \geq 0$ and $\mathcal{C}^+ \cap (-\mathcal{C}^+) = \{ 0 \}$. If $x \in X$, we write $x \geq 0$ when $x \in \mathcal{C}^+$. We define a partial order on $\mathcal{C}^+$ by letting $x \leq y$ when $y-x \in \mathcal{C}^+$.

    An operator $T : X \to X$ is said to be positive (with respect to the basis $(e_n)_{n \geq 0}$) whenever $T(\mathcal{C}^+) \subseteq \mathcal{C}^+$, that is, when $T x \geq 0$ for every $x \geq 0$.
    We write $T \geq 0$ when $T$ is a positive operator, and we write $\px$ the space of all bounded positive operators on $X$. 

    Let $(t_{i,j})_{i,j \geq 0}$ be the matrix of $T$ with respect to the basis $(e_n)_{n \geq 0}$ (with $t_{i,j} = \langle e_i^*, Te_j \rangle$ for every $i,j \geq 0$). Then $T$ is a positive operator if and only if $t_{i,j} \geq 0$ for every $i,j \geq 0$.

    We notice that if $S,T$ are two positive operators on $X$, then for every $\lambda \geq 0$, the operators $ST$, $\lambda S$ and $S+T$ are also positive. Moreover, the projections $P_N$ are positive on $X$. 

    Recall that if $M >0$, we write $\pmx$ for the set of all positive operators on $X$ with norm at most equal to $M$. In particular, if $M = 1$, the set $\ppx$ is just the positive contractions on $X$. 

    Our aim in this note is to investigate typical properties of positive contractions for the \sot \, topology and the \sote \, topology and, more precisely, we would like to know if the property to have a non-trivial invariant subspace is a typical property of positive contractions. The first step in this investigation is to check that $(\ppx, \sot)$ and $(\ppx,\sote)$ are Polish spaces. We notice that for every $M >0$, the space $\pmx$ is closed in $\bmx$ for the \sot \, topology (and hence for the \sote \, topology), and since $(\bmx,\sot)$ is Polish when $X$ is separable and $(\bmx, \sote)$ is Polish when $X^*$ is separable (see \cite[Page 256, Proposition 1.3]{Con}), we have the following result.
    
    \smallskip
    
    \bpr \label{polishsotsote}
    Let $X$ be a Banach space with a basis.
    For every $M > 0$, the space $(\pmx,\emph{\sot})$ is Polish.
    If moreover $X^*$ is separable, then the space $(\pmx,\emph{\sote})$ is Polish for every $M > 0$.
    \epr

\smallskip

\subsection{A tool for proving density results} When we study typical properties of positive contractions, we have to prove that certain sets of operators are dense in $\ppx$. The following lemma is very useful for this.

\smallskip

\blm \label{lemapprox}
Let $X$ be a Banach space with a monotone basis $(e_n)_{n \geq 0}$. Let $\mathcal{C}(X)$ be a class of operators on $X$ and define $\mathcal{C}_M(X) := \mathcal{C}(X) \cap \pmx$ for every $M > 0$.

Let $M > 0$. Suppose that the following property holds: there exists an index $p \in \Z_+$ such that for every $\varepsilon >0,$ every $N \in \Z_+$ with $N \geq p$ and every positive operator $A \in \mathcal{P}(E_N)$ with $\lVert A \rVert < M$, there exists a positive operator $T \in \mathcal{C}_M(X)$ such that
\begin{equation} \label{eq1lemapprox}
    \lVert (T-A) \, e_k \rVert < \varepsilon \quad \textrm{for every} \quad  0 \leq k \leq N.
\end{equation}

Then $\mathcal{C}_M(X)$ is dense in $(\pmx, \emph{\sot})$. If $X$ has a shrinking monotone basis $(e_n)_{n \geq 0}$ and if the condition (\ref{eq1lemapprox}) is replaced by the following condition
\begin{equation} \label{eq2lemapprox}
    \lVert (T-A) \, e_k \rVert < \varepsilon \quad \textrm{and} \quad  \lVert (T-A)^* e_k^* \rVert < \varepsilon \quad \textrm{for every} \quad 0 \leq k \leq N,
\end{equation}
then $\mathcal{C}_M(X)$ is dense in $(\pmx, \emph{\sote})$.
\elm

\smallskip

\bpf
We will prove the lemma for the \sote \,topology. 

Let $T_0 \in \pmx$, let $\varepsilon > 0$, let $x_1,...,x_s \in X$ and let $y_1^*,...,y_s^* \in X^*$. Without loss of generality, we can suppose that $\lVert T_0 \rVert < M$. We are looking for a positive operator $T \in \mathcal{C}_M(X)$ such that
\begin{equation} \label{eq1 lem1}
    \displaystyle \max_{1 \leq j \leq s} \max\{ \lVert (T-T_0)x_j \rVert, \lVert (T-T_0)^* y_j^* \rVert \} < \varepsilon.
\end{equation}
For every $1 \leq j \leq s $, there exist two indices $N_j$ and $N_j'$ such that 
$$
\left \lVert x_j - \displaystyle \sum_{k=0}^{N_j} e_k^*(x_j) e_k \right \rVert < \frac{\varepsilon}{4 M} \quad \textrm{and} \quad \; \left \lVert y_j^* - \displaystyle \sum_{k=0}^{N_j'} e_k^{**}(y_j^*) e_k^{*} \right \rVert < \frac{\varepsilon}{4 M}.
$$
Let $N_0 = \displaystyle \max_{1 \leq j \leq s}\{N_j, N_j', p \} \in \Z_+$.

We claim that if the following inequality
\begin{equation}\label{eq2 lem1}
    \max_{0 \leq k \leq N_0} \max\{ \lVert (T-T_0) e_k \rVert, \lVert (T-T_0)^* e_k^* \rVert \} < \frac{\varepsilon}{2\alpha}
\end{equation}
holds with
$$
\alpha := \max_{1 \leq j \leq s}\left\{ \displaystyle \sum_{k=0}^{N_0} \lVert x_j \rVert \lVert e_k^* \rVert, \displaystyle \sum_{k=0}^{N_0} \lVert y_j^* \rVert \lVert e_k^{**} \rVert \right\},
$$
then the inequality (\ref{eq1 lem1}) holds too. Indeed if (\ref{eq2 lem1}) holds, we have for every $1 \leq j \leq s$: 
\begin{align*}
    \lVert (T-T_0) x_j \rVert &\leq \lVert (T- T_0) (x_j - \displaystyle \sum_{k=0}^{N_j} e_k^*(x_j) e_k) \rVert + \lVert (T-T_0) \displaystyle \sum_{k=0}^{N_j} e_k^*(x_j) e_k \rVert \\
    &< \frac{\varepsilon}{2} + \displaystyle \sum_{k=0}^{N_0} \lVert e_k^*(x_j) \rVert \lVert (T-T_0) e_k \rVert \\
    &< \varepsilon
\end{align*}
and likewise we have
\begin{align*}
    \lVert(T-T_0)^* y_j^* \rVert < \varepsilon.
\end{align*}

Now for every $N \geq N_0$, we consider the positive operator $A_N := P_N T_0 P_N$. Since the basis is monotone and $\lVert T_0 \rVert < M$ we have that $\lVert A_N \rVert < M$, so there exists a positive operator $T \in \mathcal{C}_M(X)$ such that
$$
\max_{0 \leq k \leq N} \max \{ \lVert (T - A_N) e_k \rVert, \lVert (T-A_N)^* e_k^* \rVert \} < \frac{\varepsilon}{4 \alpha}.
$$
Now, for every $ 0\leq k \leq N$:
\begin{align*}
    \lVert (T- T_0) e_k \rVert &\leq \lVert (T - A_N) e_k \rVert + \lVert (P_N T_0 P_N - T_0) e_k \rVert \\
    &< \frac{\varepsilon}{4 \alpha} + \lVert (P_N - I) T_0 e_k \rVert
\end{align*}
and
\begin{align*}
    \lVert (T- T_0)^* e_k^* \rVert < \frac{\varepsilon}{4 \alpha} + \lVert (P_N-I)^* T_0^* e_k^* \rVert.
\end{align*}
Since $P_N \underset{N \to \infty}{\overset{\sote}{\longrightarrow}} I$, we can choose $N \in \Z_+$ large enough such that
$$
\max_{0 \leq k \leq N} \max \{ \lVert (P_N - I) T_0 \, e_k \rVert, \lVert (P_N-I)^* T_0^* e_k^* \rVert \} < \frac{\varepsilon}{4 \alpha},
$$
and the inequality (\ref{eq2 lem1}) follows. This proves Lemma \ref{lemapprox}.
\epf

\smallskip

\begin{remark} \label{rqapproxpos}
    The hypothesis "for every positive operator $A \in \mathcal{P}_1(E_N)$ with $\lVert A \rVert < M$" in Lemma \ref{lemapprox} can be replaced by "for every positive operator $A \in \mathcal{P}_1(E_N)$ with $\lVert A \rVert < M$ and with $\lan e_k^*, A e_l \ran > 0$ for every $0 \leq k,l \leq N$". Indeed, the operators $A_N$ in the proof of Lemma \ref{lemapprox} can be approximated in the \sote-topology by operators on $E_N$ whose matrices have positive entries.
\end{remark}

\smallskip

\begin{remark}
    Lemma \ref{lemapprox} requires the basis $(e_n)_{n \geq 0}$ to be monotone. Notice that the norm $\vertiii{.}$ defined by $\vertiii{x} = \displaystyle \sup_{N \geq 0} \lVert P_N x \rVert$ is equivalent to the norm $\lVert . \rVert$ and that when $X$ is equipped with this norm, the basis $(e_n)_{n \geq 0}$ becomes monotone.
\end{remark}

\subsection{Topological 0-1 law for positive operators} \label{ssection0-1}

    We assume in this subsection that $X = \ell_q$ with $1 \leq q < \infty$ and we consider the set 
    \begin{align*}
        \textrm{Iso}_+(X) := \{ T \in \ppx : T \; \textrm{is a surjective isometry of $X$} \}.
    \end{align*}

    It is a classical fact (see \cite[Proposition 2.f.14]{LT}) that if $1 < q \ne 2 < \infty$, every surjective isometry $T$ of $X$ has the form 
    $$T x = ( \varepsilon_n x_{\sigma(n)})_{n \geq0}, \quad \textrm{for every} \, x=(x_n)_{n \geq 0} \in X,$$ 
    where $\sigma : \Z_+ \to \Z_+$ is a bijection of $\Z_+$ and where $(\varepsilon_n)_{n \geq 0}$ is a sequence of numbers such that $\lvert \varepsilon_n \rvert = 1$ for every $n \geq 0$ . It follows that every positive surjective isometry of $X$ has the form 
    $$T x = (x_{\sigma(n)})_{n \geq0}, \quad \textrm{for every} \, x=(x_n)_{n \geq 0} \in X,$$ 
    where $\sigma : \Z_+ \to \Z_+$ is a bijection of $\Z_+$. In fact, every positive surjective isometry of $\ell_2$ also has this form. Indeed, if $T : \ell_2 \to \ell_2$ is a positive surjective isometry of $\ell_2$, then for every $i,j \geq 0 $ with $i \ne j$, we have that
    \begin{equation} \label{eqisoposl2}
        \lan T e_i, T e_j \ran = \displaystyle \sum_{k \geq 0} \lan e_k^*, T e_i \ran \, \lan e_k^*, T e_j \ran = 0.
    \end{equation}
    So using the fact that every coefficient in (\ref{eqisoposl2}) is non-negative, we have that
    \begin{equation}
        \lan e_k^*, T e_i \ran \, \lan e_k^*, T e_j \ran = 0 \quad \textrm{for every} \; k \geq 0,
    \end{equation}
    that is, $T e_i$ and $T e_j$ have disjoint supports. Now since we know that the vectors $T e_i$ for $i \geq 0$ have mutually disjoint supports, the proof given in \cite[Proposition 2.f.14]{LT} also works for the positive surjective isometries of $\ell_2$.

    From this description of the positive surjective isometries of $\ell_q$ with $1 \leq q < \infty$, it follows that Iso$_+$($X$) is a group. We say that a subset $\mathcal A$ of $\ppx$ is \emph{${\rm Iso}_+(X)\,$-$\,$invariant} if $J\mathcal AJ^{-1}=\mathcal A$ for every $J\in{\rm Iso}_+(X)$. The following result shows that every property of positive contractions that we will consider in this article is either typical or atypical when $X = \ell_q$ with $1 \leq q < \infty$.

    \smallskip

    \bpr\label{0-1} Let $X=\ell_q$ with $1\leq q<\infty$. If $\mathcal A\subseteq(\ppx,\emph{\sot})$ has the Baire property and is ${\rm Iso}_+(X)\,$-$\,$invariant, then $\mathcal A$ is either meager or comeager in $(\ppx, \emph{\sot})$. If $1 < q < \infty $ and if $\mathcal A\subseteq(\ppx,\emph{\sote})$ has the Baire property and is ${\rm Iso}_+(X)\,$-$\,$invariant, then $\mathcal A$ is either meager or comeager in $(\ppx, \emph{\sote})$
    \epr

\smallskip

    \bpf The proof given in \cite[Proposition 3.2]{GMM1}, which relies on \cite[Theorem 8.46]{Ke}, works in exactly the same way for positive contractions.
    \epf

\subsection{Similar topologies}

The aim of this subsection is to link the topologies \sot \, and \sote\, on $\ppx$ in terms of comeager sets when $X = \ell_q$ with $q > 2$. The main result of this subsection is the following.

\smallskip

\bth \label{sotsimsoteintro}
Let $X = \ell_q$ with $q > 2$. The Baire spaces $(\ppx,\emph{\sot})$ and $(\ppx, \emph{\sote})$ have the same comeager sets. 
\eth

\smallskip

The proof of Theorem \ref{sotsimsoteintro} is very similar to the proof of \cite[Theorem 3.4]{GMM2}, but we have to adapt the proof given in \cite{GMM2} to positive contractions. Indeed, the proof given in \cite{GMM2} uses \cite[Corollary 2.10]{GMM2}, which works in $\ppx$, and uses the two propositions \cite[Proposition 5.15]{GMM1} and \cite[Proposition 5.16]{GMM1}, and these two propositions can easily be adapted to positive contractions as long as \cite[Lemma 5.17]{GMM1} can be adapted to positive contractions. 

We first introduce some terminology.
We say that two topologies $\tau$ and $\tau'$ on $\mathcal{Y} = \mathcal{P}_1(\ell_q)$ are \emph{similar} if they have the same dense sets. Similar topologies have the same comeager sets (\cite[Lemma 2.1]{GMM2}), and simple examples show that the converse is not true in general (\cite[Remark 2.2]{GMM2}). We write $\mathbf{i}_{\tau, \tau'}$  the identity map from $(\mathcal{Y},\tau)$ to $(\mathcal{Y} , \tau')$ and $\mathcal{C}(\tau,\tau')$ the set of all points of continuity of this map. 

\smallskip

A vector $x \in X$ is said to be \emph{norming} for an operator $A\in \bx$ if $\Vert x\Vert=1$ and $\Vert Ax\Vert=\Vert A\Vert$. Given $N\ge 0$, we will say that an operator $A\in{\mathcal{P}}(E_{N})$ is \emph{absolutely exposing} if the set of all norming vectors for $A$ consists only of unimodular multiples of a single vector $x_{0}\in S_{E_{N}}$. We denote by ${\mathcal{E}}_1(E_N)$ the set of absolutely exposing positive operators $A\in{\mathcal P}_1(E_N)$.

The only part of the proof given in \cite{GMM2} that we have to adapt is \cite[Lemma 5.17]{GMM1}. We have to slightly modify the expression of the operators $A_\delta$ involved in the proof of \cite[Lemma 5.17]{GMM1} to obtain positive contractions. We recall that if $x$ is a vector of $X$, we denote by $\lvert x \rvert$ the positive vector of $X$ whose coordinates are the modulus of the coordinates of $x$. We thus have to prove the following lemma.

\smallskip

\begin{lemma}\label{lem517}
The set ${\mathcal{E}}_1(\en)$ is dense in $\mathcal{P}_{1}(\en)$.
 \end{lemma}

\smallskip

\bpf
Let $A\in\mathcal{P}_{1}(\en)$ with $A\neq 0$ and $\Vert A\Vert<1$. Let $x_{0}\in \en$ be such that $\lVert x_0 \rVert = 1$ and $\Vert Ax_{0}\Vert=\Vert A\Vert$. Since the vector $\lvert x_0 \rvert$ is also a norming vector for $A$, we can suppose without loss of generality that $x_0 \geq 0$. 
 
 By the Hahn-Banach theorem, there exists a functional $x_{0}^{*}\in\en^*$  such that $\Vert x_{0}^*\Vert=\pss{x_{0}^{*}}{x_{0}}=1$. We can also suppose that $x_0^*$ is a positive functional. Indeed if $x_0^* = \displaystyle \sum_{j = 0}^N \beta_j e_j^*$, we consider the positive functional defined by $y_0^* = \displaystyle \sum_{j=0}^N \lvert \beta_j \rvert e_j^*$. Using Hölder's inequality, we can show that $\lVert y_0^* \rVert \leq 1$, and moreover, we have that $ \lvert \lan x_0^*, x_0 \ran \rvert \leq \lan y_0^* , x_0 \ran \leq \lVert y_0^* \rVert \leq 1$, so $\lVert y_0^* \rVert = 1 = \lan y_0^* , x_0 \ran$.
 
 Let $R_{0}$ be the positive rank 1 operator on $\en$ defined by $R_{0}(x):=\pss{x_{0}^{*}}{x}\,Ax_{0}$, for every $x\in\en$;  and for any $\delta >0$, let 
 $A_{\delta }$ be the operator defined by 
 $A_{\delta }:=A+\delta R_{0}$. The operators $A_\delta$ are now positive. As in the proof of \cite[Lemma 5.17]{GMM1}, we can prove that $A_{\delta }$ is absolutely exposing and given $\varepsilon >0$, one can choose $\delta >0$ so small that $\Vert A_{\delta }\Vert<1$ and $\Vert A-A_{\delta }\Vert<\varepsilon $, because $\lVert A \rVert < 1$. This proves that ${\mathcal{E}}_1(\en)$ is dense in $\mathcal{P}_{1}(\en)$.
\epf

Since the proof of \cite[Theorem 3.4]{GMM2} can now be adapted to positive contractions thanks to Lemma \ref{lem517}, we can state the following results.

\smallskip

\bth\label{sotsim} Let $X = \ell_q$ with $q > 2$. Then the topologies \emph{\sot}  and \emph{\sote} are similar on $\mathcal{P}_1(X)$. 
\eth

\smallskip

\bco \label{cortopsimsot}
Let $X = \ell_q$ with $q >2$. The Baire spaces $(\ppx, \emph{\sot})$ and $(\ppx, \emph{\sote})$ have the same comeager sets.
\eco

\section{Some elementary properties of typical positive contractions}\label{Section3}

In this section, we study some elementary properties of a typical positive contraction on $\ell_q$ for the \sot \, and the \sote \,topologies. These properties will be useful in the next section to prove that a typical positive contraction on $\ell_1$ and on $\ell_2$ has a non-trivial invariant subspace. 

\smallskip

\bpr \label{normspec} Let $X = \ell_q$ with $1 \leq q < \infty$. A typical $T \in (\ppx, \emph{\sot})$ has the following properties:
\begin{enumerate}[\rm (a)]
    \item $\lVert T \rVert = 1$;
    \item $\Vert T^nx\Vert\to 0$ as $n\to\infty$ for all $x\in X$;
    \item $T$ is not invertible;
    \item $\sigma_{ap}(T) = \sigma(T) = \overline{\D}$.
\end{enumerate}
If $1 < q < \infty$, a typical $T \in (\ppx, \emph{\sote})$ also has these properties. 
\epr

\smallskip

\bpf
The property (a) follows from the fact that the set 
$$\mathcal{A} := \{ T \in \ppx : \lVert T \rVert = 1\}$$ 
can be written as 
$$ 
\mathcal{A} = \displaystyle \bigcap_{k \geq 1} \bigcup_{x \in S_X} \{ T \in \ppx : \lVert Tx \rVert > 1 - \frac{1}{k}\}.
$$
So $\mathcal{A} $ is a \sot-$G_\delta$ subset of $\ppx$ (and hence a \sote-$G_\delta$ subset of $\ppx$). 

\noindent Moreover, if $T$ is a positive contraction and if we set $T_N = P_N T P_N + Q_N$ for every $N \geq 0$, then $T_N$ belongs to $\mathcal{A}$ and we easily see that $T_N \underset{n \to \infty}{\overset{\sote}{\longrightarrow}}T$, so $\mathcal{A}$ is dense in $\ppx$ for both the \sot \, and the \sote\, topologies.

The proof of properties (b), (c) and (d) can easily be adapted from \cite[Propositions~3.7 to 3.9]{GMM1} to positive contractions.
\epf

Recall that an operator $T \in \bx$ is said to be hypercyclic if there is a vector $x \in X$ such that the orbit $O(x,T) := \{ T^n x : n \in \N \}$ is dense in $X$. Such a vector is said to be hypercyclic for $T$. We refer to \cite{BM} and \cite{GEP} for background on hypercyclicity.

\smallskip

\bpr \label{prophypercycl} Let $X = \ell_q$ with $1 \leq q < \infty$. For any $M > 1$, the set 
$$\{ T \in \pmx : T \;\textrm{is hypercyclic} \}$$ 
is comeager in $(\pmx, \emph{\sot})$ and in $(\pmx, \emph{\sote})$.
In particular, a typical positive contraction $T \in (\ppx,\emph{\sot})$ is such that $2T$ is hypercyclic. 

If $1 < q < \infty$, a typical $T \in (\ppx,\emph{\sote})$ is such that $2T$ and $(2T)^*$ are hypercyclic.

\epr

\smallskip

\bpf 
Using \cite[Theorem 5.41]{BM} and replacing the property of being mixing in the Gaussian sense by the property of being hypercyclic in \cite[Lemma 2.8]{GMM}, we notice that the proof given in \cite[Proposition 2.3]{GMM} works for any $1 \leq q < \infty$ in $\mathcal{P}_M(X)$. When $1 < q < \infty$, the map $T \mapsto T^*$ is a homeomorphism from $(\mathcal{P}_2(X), \sote)$ to $(\mathcal{P}_2(X), \sote)$, so an \sote-typical $T \in \ppx$ is such that $(2T)^*$ is hypercyclic too. 
\epf

\smallskip
Using the fact that a hypercyclic operator $T$ is such that $\sigma_p(T^*) = \emptyset$, we have the following result.

\smallskip

\bco \label{corvalprop}
Let $X = \ell_q$ with $1 \leq q < \infty$. An \emph{\sot}-typical $T \in \ppx$ is such that $T^*$ has no eigenvalue. If $1 < q < \infty$, an \emph{\sote}-typical $T \in \ppx$ is such that $T$ and $T^*$ have no eigenvalue.
\eco

\smallskip
Using Corollary \ref{cortopsimsot}, we also have the following result.

\smallskip

\bco \label{corvalpropbis}
Let $X = \ell_q$ with $q > 2$. An \emph{\sot}-typical $T \in \ppx$ is such that $T$ has no eigenvalue.
\eco

\smallskip

The next corollary is an exact analogue of \cite[Proposition 3.9]{GMM1}.

\smallskip
\bco 
If $X = \ell_q$ with $1 \leq q < \infty$, a typical $T \in (\ppx, \emph{\sot})$ is such that $T-\lambda$ has dense range for every $\lambda \in \C$. If $q >1$, a typical $T \in (\ppx, \emph{\sote})$ is such that $T-\lambda$ has dense range for every $\lambda \in \C$.
\eco

\smallskip
\bpf
This is clear from Proposition \ref{prophypercycl} because a hypercyclic operator $T$ is such that $T-\lambda$ has dense range for every $\lambda \in \C$.
\epf

Our next step is to investigate whether a typical $T \in \ppx$ is such that $T^*$ is an isometry or not. The following fact that we already proved in Subsection \ref{ssection0-1} will be useful for this.

\smallskip

\begin{fact} \label{fact1iso}
    If $T$ is a positive isometry of $X = \ell_q$ with $1 \leq q  < \infty$, then the vectors $Te_i$ for $i \geq 0$ have mutually disjoint supports.
\end{fact}

\smallskip

By \cite[Proposition 5.15]{EM} and \cite[Corollary 3.5]{GMM1}, we know that an \sot-typical $T \in \mathcal{B}_1(\ell_2)$ is such that $T^*$ is an isometry. We prove that this is no longer the case for a typical positive contraction on $\ell_2$. 

\smallskip

\bpr \label{propcoisolq}
Let $X = \ell_q$ with $1 < q <\infty$. A typical $T \in (\ppx,\emph{\sot})$ (resp. $T \in (\ppx,\emph{\sote})$) is such that $T^*$ is not an isometry. 
\epr

\smallskip

\bpf
The arguments given in \cite[Proposition 5.1]{GMM1} also apply in this case. We will detail a bit the denseness argument. Let 
$$\mathcal{I}^* := \{ T \in \ppx : T^* \; \textrm{is an isometry}\}.$$ 

By Fact \ref{fact1iso}, we have that

$$\mathcal{A} :=\displaystyle \bigcup_{j \geq 0} \{ T \in \ppx : \langle e_0^*, T e_j \rangle \ne 0 \;\textrm{and} \; \langle e_1^* , T e_j \rangle \ne 0  \} \subseteq \ppx \setminus \mathcal{I}^*$$
and the set $\mathcal{A}$ is \sot-open in $\ppx$, so $\mathcal{A}$ is a \sot-$G_\delta$ of $\ppx$ and hence \sote-$G_\delta$. It remains to prove that $\mathcal{A}$ is dense in $(\ppx, \sote)$.
Let $\varepsilon > 0$, let $T \in \ppx$ with $\lVert T \rVert < 1$ and let $x_1,...,x_n \in X$ and $y_1^*,..., y_n^* \in X^*$. We have to find a positive contraction $S$ in the set $\mathcal{A}$ such that
\begin{equation} \label{eqensA}
    \displaystyle \max_{1 \leq l \leq n} \max\{ \lVert (T-S) \,x_l \rVert, \lVert (T-S)^* y_l^* \rVert \} < \varepsilon.
\end{equation}
Consider the positive operator $S_\delta$ defined by $S_\delta(x) = T x + \delta \, \lan e_0^*, x \ran \, (e_0 + e_1)$ for every $x \in X$, where $\delta$ is a positive number that we will define later on. We have that 
$$
\lan e_0^*, S_\delta \, e_0 \ran \, \lan e_1^*, S_\delta \, e_0 \ran \geq \delta^2 > 0.
$$
For every $x \in X$, we have that
$$
\lVert S_\delta \, x \rVert \leq (\lVert T \rVert + 2 \delta) \lVert x \rVert
$$
and for every $1 \leq l \leq n$, we have that
$$
\lVert (T - S_\delta) \, x_l \rVert \leq 2 \delta \lVert x_l \rVert \quad \textrm{and} \quad \lVert (T - S_\delta)^* y_l^* \rVert \leq 2 \delta \lVert y_l^* \rVert.
$$
If we choose $\delta > 0$ such that
$$
\delta < \frac{1- \lVert T \rVert}{2}, \quad 2 \delta  \lVert x_l \rVert < \varepsilon \quad \textrm{and} \quad 2 \delta  \lVert y_l^* \rVert < \varepsilon \quad \textrm{for every} \; 1 \leq l \leq n,
$$
then the operator $S_\delta$ is a positive contraction of $\mathcal{A}$ satisfying (\ref{eqensA}). This concludes the proof of Proposition \ref{propcoisolq}. 
\epf

The case $X = \ell_1$ is very different from the case $X = \ell_q$ with $q >1$. We have the following result which will be useful in Section \ref{Section4} to identify the point spectrum of an \sot-typical positive contraction on $\ell_1$.

\smallskip

\bpr \label{propsurjl1}
Let $X = \ell_1$. An \emph{\sot}-typical $T \in \ppx$ is such that $T^*$ is a non-surjective isometry and such that $T-\lambda$ is surjective for every $\lambda \in \D$.
\epr

\bpf
The proof works exactly as in \cite[Theorem 4.1]{GMM1}. Indeed, the set
$$
\mathcal{I}^* := \{ T \in \ppx : T^* \; \textrm{is an isometry} \}
$$
is a \sot-$G_\delta$ of $\ppx$ (see \cite[Theorem 4.1]{GMM1}) and the set $\mathcal{I}^*$ is also \sot-dense in $\ppx$ since all the operators involved in the proof of \cite[Theorem 4.1]{GMM1} are positive. The second part of the proof immediately follows as in \cite[Theorem 4.1]{GMM1}.
\epf

Finally, we describe the essential spectrum of a typical positive contraction on $X = \ell_q$ with $1 \leq q < \infty$. This will be useful to prove that a typical positive contraction on $\ell_1$ has a non-trivial invariant subspace. Recall that an operator $T \in \bx$ is Fredholm if its kernel is finite-dimensional and its range has finite codimension, and it is upper semi-Fredholm if its range is closed and its kernel is finite-dimensional. Every Fredholm operator is upper semi-Fredholm. An operator is semi-Fredholm if it is upper semi-Fredholm or if its range has finite codimension. We denote by $\phip$ the set of all upper semi-Fredholm operators on $X$. The set $\phip$ is norm-open in $\bx$. We refer to \cite{Mul} for background on Fredholm operators. The first lemma that we will use is the following analogue of \cite[Fact 7.14]{GMM2}.

\smallskip

\blm \label{lemfred1}
Let $X = \ell_q$ with $1 \leq q < \infty$ and let $\lambda \in \overline{\D}$. A typical $T \in (\ppx,\emph{\sot})$ (resp.  $T \in (\ppx,\emph{\sote})$ when $1 < q < \infty$)  has the following property:

For every $\varepsilon >0$ and every $n \geq 1$, there exists a subspace $E$ of $X$ with $n < \textrm{dim}(E) < \infty$ such that $\lVert \left. (T-\lambda) \right|_E \rVert < \varepsilon$.  
\elm

\smallskip

\bpf
Let $\mathcal{G}$ be the set of all operators $T\in \ppx$ satisfying this property. Then 
$$ 
\displaystyle \mathcal{G} = \bigcap_{\substack{p \geq 0 \\ n \geq 1}} \{ T\in \ppx : \exists \, E \; \textrm{subspace of $X$}, \; n < \textrm{dim}(E) < \infty, \; \lVert \left. (T-\lambda) \right|_E \rVert < 2^{-p} \}.
$$
For every subspace $E$ of $X$ satisfying $n < \textrm{dim}(E) < \infty$, the set 

$$\Lambda_{E,p} :=\{ T \in \ppx : \lVert \left. (T-\lambda) \right|_E \rVert < 2^{-p} \} $$ 
is \sot-open (see \cite[Fact 7.14]{GMM2}), so $\mathcal{G}$ is \sot-$G_\delta$ in $\ppx$ and hence \sote-$G_\delta$ in $\ppx$. 

Moreover, if $T$ is a positive contraction and if we set $T_N = P_N T P_N + \lambda \, Q_N$ for every $N \geq 0$, then $T_N$ is a positive contraction which belongs to $\mathcal{G}$, and $T_N \underset{N \to \infty}{\overset{\sote}{\longrightarrow}}T$. So $\mathcal{G}$ is dense in $\ppx$ for the \sot\, and the \sote \,topologies.
 \epf

 \smallskip

With Lemma \ref{lemfred1}, we obtain the following description of the essential spectrum of a typical positive contraction. 

\smallskip

\bpr  \label{propFredh}
Let $X = \ell_q$ with $1 \leq q < \infty$. A typical $T \in (\ppx,\emph{\sot})$ (resp. $T \in (\ppx,\emph{\sote})$ when $1 < q < \infty$) is such that $T-\lambda$ is not upper semi-Fredholm for every $\lambda \in \overline{\D}$ and such that $\sigma_{ess}(T) = \overline{\D}$. 
\epr

\smallskip

\bpf
The proof given in \cite[Proposition 7.13]{GMM2} works in our case since Lemma \ref{lemfred1} is working in $\ppx$ for the topologies \sot\, and \sote.
\epf

\section{Invariant subspaces and typicality}\label{Section4}

We now come to our main goal in this paper, which is to investigate whether a typical $T \in (\ppx, \sot)$ (resp. $T \in (\ppx, \sote)$) has a non-trivial invariant subspace. We will in fact see that this is not an easy question to answer. A first observation is that the property of having a non-trivial invariant subspace is either typical or atypical (\cite[Corollary~3.3]{GMM1}).

\smallskip

\bpr Let $X=\ell_q$ with $1 \leq q < \infty $. Either a typical $T\in (\ppx,\emph{\sot})$ (resp. $T \in (\ppx, \emph{\sote})$ when $1 < q < \infty$) has a non-trivial invariant subspace, or a typical $T\in (\ppx,\emph{\sot})$ (resp. $T \in (\ppx, \emph{\sote})$ when $1 < q < \infty$) does not have a non-trivial invariant subspace.
\epr

\smallskip

An important result from Brown, Chevreau and Pearcy (\cite{BCP2}) states that every contraction on a Hilbert space whose spectrum contains the unit circle has a non-trivial invariant subspace. Since a typical $T \in (\mathcal{P}_1(\ell_2),\sot)$ (resp. $T \in (\mathcal{P}_1(\ell_2),\sote)$) is such that $\sigma(T) = \overline{\D}$ by Proposition \ref{normspec}, we have the following result as in \cite[Corollary 7.3]{GMM1}.

\smallskip

\bth
A typical $T \in (\mathcal{P}_1(\ell_2), \emph{\sot})$ (resp. $T \in(\mathcal{P}_1(\ell_2), \emph{\sote})$) has a non-trivial invariant subspace.
\eth

\subsection{Invariant subspace of typical positive contractions on $\ell_1$}

In this subsection, we consider the case where $X = \ell_1$. We will prove that a typical positive contraction on $X = \ell_1$ has a non-trivial invariant subspace, and even has eigenvalues. 

By Proposition \ref{normspec} and Proposition \ref{propsurjl1}, we know that a typical $T \in (\ppx,\sot)$ is such that $T-\lambda$ is surjective for every $\lambda \in \D$ and such that $\sigma(T) = \overline{\D}$. This implies that a typical $T \in (\ppx,\sot)$ is such that $T-\lambda$ is not injective for every $\lambda \in \D$. Thus, we have the following results.

\smallskip

\bth \label{thvalpropl1}
Let $X = \ell_1$. A typical $T \in (\ppx,\emph{\sot})$ is such that $\sigma_p(T) = \D$ and such that \emph{dim(Ker($T-\lambda$))} $= \infty$ for every $\lambda \in \D$.
\eth

\smallskip

\bpf[Proof of Theorem \ref{thvalpropl1}]

The proof is motivated by \cite[Remark 4.5]{GMM1}.  
By Propositions \ref{propsurjl1} and \ref{propFredh}, a typical $T \in (\ppx,\sot)$ is not Fredholm and is surjective, so a typical $T \in (\ppx,\sot)$ is such that dim(Ker($T$)) $= \infty$. But a typical $T \in (\ppx,\sot)$ is such that $T-\lambda$ is semi-Fredholm for every $\lambda \in \D$ (because it is surjective). By the continuity of the Fredholm index, Ind($T-\lambda$) does not depend on $\lambda \in \D$. So a typical $T \in (\ppx,\sot)$ is such that Ind($T-\lambda$) = dim(Ker($T-\lambda$)) $= \infty$ for every $\lambda \in \D$.
\epf

\smallskip

\bco \label{corinvsubl1}
Let $X = \ell_1$. A typical $T \in (\ppx, \emph{\sot})$ has a non-trivial invariant subspace.
\eco

\smallskip

Hence, the cases $X = \ell_2$ and $X = \ell_1$ are fully understood. We will now see that the problem is more difficult in the other cases, which is not surprising because the problem is still open in $\bbx$ when $X = \ell_q$ with $1 < q \ne 2 < \infty$ for the topologies $\sot$ and $\sote$.

\subsection{Invariant subspace of typical positive contractions on a Banach space with a basis.}

We now focus on the case where $X$ is a Banach space with a basis. In particular, it includes the case $X = \ell_q$ with $1 < q \ne 2 < \infty$.

We start with a similar result to \cite[Proposition 5.24]{GMM1} in the case where $X = \ell_q$ with $1 \leq q < \infty$. 

Recall that an operator $T \in \bx$ is polynomially bounded if there exists $C > 0$ such that for every complex polynomial $P$:
$$
\lVert P(T) \rVert \leq C \sup_{\lvert z \rvert = 1} \lvert P(z) \rvert.
$$
Every contraction on $\ell_2$ is polynomially bounded by Von Neumann's inequality. A remarkable result of Ambrozie and Müller \cite{AM} states that every polynomially bounded operator $T \in \bx$ such that $T^n \underset{n \to \infty}{\overset{\sot}{\longrightarrow}}0$ and $\T \subseteq \sigma(T)$ has a non-trivial invariant subspace. We know that a typical $T \in \ppx$ is such that $T^n \underset{n \to \infty}{\overset{\sot}{\longrightarrow}}0$ and $\sigma(T) = \overline{\D}$ by Proposition \ref{normspec}, so it is natural to ask whether a typical $T \in \ppx$ is polynomially bounded or not.

\smallskip

\bpr
Let $X = \ell_q$ with $1 \leq q \ne 2 < \infty$. A typical $T\in (\ppx,\emph{\sot})$ (resp. $T\in (\ppx,\emph{\sote})$ when $1 < q \ne 2 < \infty$ ) is not polynomially bounded. 
\epr

\smallskip

\bpf
The proof given in \cite[Proposition 5.24]{GMM1} can easily be adapted to positive contractions on $X$.
\epf

Recall that a positive operator on a Banach space with a basis is said to be satisfying the AAB criterion if it satisfies the following theorem (already stated in the introduction) due to Abramovich, Aliprantis and Burkinshaw (\cite[Theorem 2.2]{AAB}).

\smallskip

\bth \label{thinvsubpospart4}
Let $X$ be a Banach space with a basis and let $T$ be a positive operator on $X$. If there exists a non-zero positive operator $A$ on $X$ which is quasinilpotent at a certain non-zero positive vector of $X$ and such that $A T = T A$, then $T$ has a non-trivial invariant subspace. 
\eth

\smallskip

Recall also that by \cite[Theorem 2.2]{AAB2}, any positive operator on $\ell_q$ satisfying the AAB criterion has a non-trivial closed invariant ideal, that is, a closed vector subspace $V$ of $\ell_q$ such that $\lvert x \rvert \leq \lvert y \rvert$ and $y \in V$ imply $x \in V$, for every $x,y \in \ell_q$. Moreover, by \cite[Proposition 1.2]{RT}, a positive operator on a Banach space with an unconditional basis $(e_n)_{n \geq 0}$ has no non-trivial closed invariant ideals if and only if it satisfies the following condition:
$$
\forall \, i \ne j \in \Z_+, \, \exists n \in \Z_+ : \, \lan e_j^*, T^n e_i \ran > 0.
$$

\smallskip

The following result shows that a typical positive contraction on a Banach space with an unconditional basis has no non-trivial closed invariant ideals.

\smallskip

\bpr \label{propidealinvsec4}
Let $X$ be a Banach space with an unconditional basis $(e_n)_{n \geq 0}$. A typical $T \in (\ppx, \emph{\sot})$ (resp. $T \in (\ppx, \emph{\sote})$ when $X^*$ is separable) has no non-trivial closed invariant ideals.  
\epr

\bpf
Consider the set
$$
\mathcal{G} := \{ T \in \ppx : \forall \, i \ne j \in \Z_+, \exists \, n \in \Z_+ \; \textrm{such that} \; \lan e_j^*, T^n e_i \ran > 0 \}.
$$
Then
$$
\mathcal{G} = \displaystyle \bigcap_{\substack{i,j \geq 0 \\ i \ne j}} \bigcup_{n \geq 0} \{ T \in \ppx : \, \lan e_j^*, T^n e_i \ran > 0\} = \displaystyle \bigcap_{\substack{i,j \geq 0 \\ i \ne j}} \mathcal{G}_{i,j} \,,
$$
with
$$\mathcal{G}_{i,j} := \displaystyle \bigcup_{n \geq 0} \{ T \in \ppx : \lan e_j^*, T^n e_i \ran > 0 \} \quad \textrm{for every} \; i,j \geq 0 \; \textrm{with} \; i \ne j. $$
Each $\mathcal{G}_{i,j}$ is easily seen to be \sot-open (and hence \sote-open) because the map $T \mapsto T^n$ is continuous on $\ppx$ for the \sot \,topology. 
Let us now show that every $\mathcal{G}_{i,j}$ is \sote-dense in $\ppx$. Let $\varepsilon > 0$, let $T \in \ppx$ with $\lVert T \rVert < 1$ and let $x_1,...,x_n \in X$ and $y_1^*,..., y_n^* \in X^*$. We have to find a positive contraction $S$ in the set $\mathcal{G}_{i,j}$ such that
\begin{equation} \label{eqidealinvsec4}
    \displaystyle \max_{1 \leq l \leq n} \max\{ \lVert (T-S) \,x_l \rVert, \lVert (T-S)^* y_l^* \rVert \} < \varepsilon.
\end{equation}
Consider the positive operator $S_\delta$ defined by $S_\delta(x) = T x + \delta \, \lan e_i^*, x \ran \, e_j$ for every $x \in X$, where $\delta$ is a positive number that we will define later on. We have that 
$$
\lan e_j^*, S_\delta \, e_i \ran \geq \delta > 0.
$$
For every $x \in X$, we have that
$$
\lVert S_\delta \, x \rVert \leq (\lVert T \rVert + \delta \lVert e_i^* \rVert \lVert e_j \rVert) \lVert x \rVert
$$
and for every $1 \leq l \leq n$, we have that
$$
\lVert (T - S_\delta) \, x_l \rVert \leq \delta \lVert x_l \rVert \lVert e_i^* \rVert \lVert e_j \rVert \quad \textrm{and} \quad \lVert (T - S_\delta)^* y_l^* \rVert \leq \delta \lVert y_l^* \rVert \lVert e_i^* \rVert \lVert e_j \rVert.
$$
If we choose $\delta > 0$ such that
$$
\delta < \frac{1- \lVert T \rVert}{\lVert e_i^* \rVert \lVert e_j \rVert}, \quad \delta  \lVert x_l \rVert \lVert e_i^* \rVert \lVert e_j \rVert < \varepsilon \quad \textrm{and} \quad \delta  \lVert y_l^* \rVert \lVert e_i^* \rVert \lVert e_j \rVert < \varepsilon \quad \textrm{for every} \; 1 \leq l \leq n,
$$
then the operator $S_\delta$ is a positive contraction of $\mathcal{G}_{i,j}$ satisfying (\ref{eqidealinvsec4}). This proves that each $\mathcal{G}_{i,j}$ is $\sote$-dense in $\ppx$ and this concludes the proof of Proposition \ref{propidealinvsec4}. 
\epf

\smallskip

In particular, we obtain the following result in the case where $X = \ell_q$.

\smallskip

\bco \label{coraabSec4}
Let $X = \ell_q$ with $1 \leq q < \infty$. A typical $T \in (\ppx, \emph{\sot})$ (resp. $T \in (\ppx, \emph{\sote})$ when $1 < q < \infty$) does not satisfy the AAB criterion.
\eco

\smallskip

Let $X$ be a Banach space with a basis. Since any positive operator on $X$ that is quasinilpotent at a non-zero positive vector of $X$ satisfies the AAB criterion, we immediately get the following result. 

\smallskip

\bco \label{corqnubcondsec4}
Let $X = \ell_q$ with $1 \leq q < \infty$.  A typical $T \in (\ppx, \emph{\sot})$ (resp. $T \in (\ppx, \emph{\sote})$ when $1 < q < \infty$) is not quasinilpotent at any non-zero positive vector of $X$. 
\eco

\smallskip

We now would like to extend Corollary \ref{coraabSec4} to other Banach spaces with a basis. The following lemma will be very useful for what follows.

\smallskip

\blm \label{lemqn}
Let $X$ be a Banach space with a basis $(e_n)_{n \geq 0}$, let $A$ be a non-zero positive operator on $X$ and let $y \in X$ be such that $y \geq 0 $ and $y \ne 0$. Let $j \in \Z_+$ be such that $\langle e_j^*, y \rangle := \alpha_j > 0$. If 
$$a_{j,j} := \langle e_j^*, A e_j \rangle > 0 ,$$
then the operator $A$ is not quasinilpotent at $y$.
\elm

\smallskip

\bpf

First, we remark that $A e_j \geq a_{j,j} \, e_j$. Since  
$$
\lan e_j^*, A y \ran = \displaystyle \sum_{k \geq 0} \lan e_j^*, Ae_k \ran \, \lan e_k^*, y \ran \,,
$$
we also have $\langle e_j^* , Ay  \rangle  \geq a_{j,j} \,\alpha_j$.
It follows that $Ay \geq \alpha_j \, a_{j,j} \, e_j$ and by induction, we easily get that
\begin{align*}
    A^k y &\geq \alpha_j \,  a_{j,j}^k \,  e_j \quad \textrm{for every} \; k \geq 1.
    \end{align*}

In particular, we get that
\begin{align*}
    \lVert A^k y \rVert \lVert e_j^* \rVert \geq \lan e_j^*, A^k y \ran \geq \alpha_j \, a_{j,j}^k \quad \textrm{for every} \; k \geq 1
\end{align*}
and thus
$$
\liminf_{k \to \infty} \lVert A^k y \rVert^{1/k} \geq a_{j,j} > 0.
$$
This shows that the operator $A$ cannot be quasinilpotent at $y$ and this concludes the proof of Lemma \ref{lemqn}. 
\epf

The main result of this paper is the following generalization of Corollary \ref{coraabSec4}.

\smallskip

\bth \label{thaabreflexivesec4}
Let $(X, \lVert . \rVert)$ be a reflexive Banach space with a monotone basis $(e_n)_{n \geq 0}$. A typical $T \in ( \ppx, \emph{\sote} )$ does not satisfy the AAB criterion.
\eth

\smallskip

\bpf
We denote by $\mathcal{F}$ the set 
$$ 
\mathcal{F} = \{ T \in \ppx : T \; \textrm{satisties the AAB criterion}\}.
$$

By Lemma \ref{lemqn}, we have
\begin{align*}
    \mathcal{F} &\subseteq \bigcup_{p \geq 0} \{ T \in \ppx : \exists A \in \ppx \; \textrm{such that} \; A \ne 0, \; A T = T A \; \textrm{and} \; \langle e_p^*, A e_p \rangle = 0 \} \\
    &\subseteq \bigcup_{p \geq 0} \bigcup_{i,j \geq 0} \bigcup_{\eta \in \Q_{>0}} \{ T \in \ppx : \exists A \in \ppx \; \textrm{such that} \; \lan e_j^*, A e_i \ran \geq \eta, \; A T = T A \; \textrm{and} \; \langle e_p^*, A e_p \rangle = 0 \} \\
    &\subseteq \bigcup_{p \geq 0} \bigcup_{i,j \geq 0} \bigcup_{\eta\in \Q_{>0}} \mathcal{F}_{i,j,\eta,p}
\end{align*}
with
$$
\mathcal{F}_{i,j,\eta,p} := \{ T \in \ppx : \exists A \in \ppx \; \textrm{such that} \; \lan e_j^*, A e_i \ran \geq \eta, \; A T = T A \; \textrm{and} \; \langle e_p^*, A e_p \rangle = 0 \}.
$$
\smallskip
We first prove the following fact.
\smallskip

\begin{fact} \label{factuno}
For each $i,j,p \geq 0$ and $\eta \in \Q_{>0}$, the set $\mathcal{F}_{i,j,\eta,p}$ is \sote-closed in $\ppx$.
\end{fact}

\begin{proof}[Proof of Fact \ref{factuno}]
Let $(T_k)_{k \geq0 } \subseteq \mathcal{F}_{i,j,\eta,p}$ be such that $T_k \underset{k \to \infty}{\overset{\sote}{\longrightarrow}} T$ with $T \in \ppx$. For every $k \geq 0$, there exists an operator $A_k \in \ppx$ such that $A_k T_k = T_k A_k, \; \lan e_j^*, A_k e_i  \ran \geq \eta$ and $\langle e_p^*, A_k e_p \rangle = 0$. Since $\ppx$ is \wot-compact (see \cite[Page 275, Proposition 5.5]{Con}), we can suppose without loss of generality that there exists an operator $A \in \ppx$ such that $A_k \underset{k\to \infty}{\overset{\wot}{\longrightarrow}}A$. We immediately get that $\langle e_p^* , A e_p \rangle = 0$ and that $\lan e_j^*, A e_i \ran \geq \eta$. Let us now show that $A T = T A$.

Let $y^* \in X^*$ and $x \in X$. We have 
\begin{align*}
\langle y^*, A_k T_k x \rangle = \langle y^*, A_k (T_k - T) x \rangle + \langle y^* , A_k Tx \rangle \quad \textrm{for every} \; k \geq 0.    
\end{align*}
Since the sequence $(A_k)_{k \geq 0}$ is bounded and since $T_k \underset{k \to \infty}{\overset{\sot}{\longrightarrow}} T $, we have that $$\lan y^*, A_k (T_k - T) x \ran \underset{k \to \infty}{\longrightarrow} 0,$$ and since $A_k \underset{k \to \infty}{\overset{\wot}{\longrightarrow}} A $, we have that $$\lan y^*, A_k T x \ran \underset{k \to \infty}{\longrightarrow} \lan y^*, A T x \ran.$$
It follows that $\langle y^*, A_k T_k x \rangle \underset{k\to \infty}{\longrightarrow} \langle y^*, A T x \rangle$.

Likewise, we have
\begin{align*}
    \langle y^*, T_k A_k x \rangle = \langle (T_k^* - T^*) y^* , A_k x \rangle + \langle T^* y^* , A_k x \rangle \quad \textrm{for every } \; k \geq 0,
\end{align*}
so using the fact that $T_k^* \underset{k \to \infty}{\overset{\sot}{\longrightarrow}} T^*$, that $A_k \underset{k \to \infty}{\overset{\wot}{\longrightarrow}} A $ and that the sequence $(A_k)_{k \geq 0}$ is bounded, we can prove by the same way that $\langle y^*, T_k A_k x \rangle \underset{k\to \infty}{\longrightarrow} \langle y^*, T A x \rangle$. It follows that
$$
\langle y^* , AT x \rangle = \langle y^*, TA x \rangle
$$
for every $y^* \in X^*$ and $x \in X$, and thus $A T = T A$. This proves that $T \in \mathcal{F}_{i,j,\eta,p}$ and concludes the proof of Fact \ref{factuno}.
\end{proof}

Let us now prove the following proposition.

\smallskip

\begin{proposition} \label{propdens}
    For each $i,j,p \geq 0$ and $\eta \in \Q_{>0}$, the set $\mathcal{F}_{i,j,\eta,p}$ has empty interior in $\ppx$ for the \emph{\sote} topology.
\end{proposition}

\smallskip
\begin{proof}[Proof of Proposition \ref{propdens}]
Using Lemma \ref{lemapprox}, we will prove that the set $\ppx \setminus \mathcal{F}_{i,j,\eta,p}$ is dense in $\ppx$ for the \sote \,topology.

Let us fix $i,j,p \in \Z_+$ and $\eta \in \Q_{>0}$.

    Let $N \in \Z_+$ be such that $N \geq p$ and let $M \in \mathcal{P}_1(E_N)$ with $\lVert M \rVert < 1$. Let $\varepsilon > 0$. Let also $u = e_0 + ... + e_{p+N+1}$ and let us write $P : X \to E_N$ the canonical projection instead of $P_N$. By Remark \ref{rqapproxpos}, we can suppose that $\lan e_k^*, M e_l \ran > 0$ for all indices $0 \leq k,l \leq N$. 
    
    Let $T$ be the positive operator defined as follows: 
    $$ 
    Tx = M P x + \delta \, \lan e_{N+p+1}^* , x \ran \, u + Sx \quad \textrm{for every} \; x \in X.
    $$
    In this expression, $S$ is the operator on $X$ defined by
    $$
    S \left(\displaystyle \sum_{k=0}^{\infty} e_k^*(x) e_k \right) = \displaystyle \sum_{k=0}^{N+p+1} \delta \, e_k^*(x) \, e_{k+N+1} + \displaystyle \sum_{k > N+p+1} \delta_{k-N-p-1} \, e_k^*(x) \, e_{k+N+1} \,,
    $$
    where $(\delta_k)_{k \geq 1}$ is a sequence of positive real numbers satisfying
    $$
    0 < \displaystyle \sum_{k > N+p+1} \delta_{k-N-p-1} \lVert e_k^* \rVert \lVert e_{k+N+1} \rVert < 1 - \lVert M \rVert \,,
    $$
    where $\delta > 0 $ is such that 
$$
0 < \delta <  \frac{1- \lVert M \rVert - \displaystyle \sum_{k > N+p+1} \delta_{k-N-p-1} \lVert e_k^* \rVert \lVert e_{k+N+1} \rVert}{\lVert u \rVert \lVert e_{N+p+1}^* \rVert + \displaystyle \sum_{k = 0}^{N+p+1} \lVert e_k^* \rVert \lVert e_{k+N+1} \rVert} \, ,
$$
and where
$$\delta \lVert e_{N+1+k} \rVert < \varepsilon, \quad \delta \lVert e_{N+1+p} ^* \rVert < \varepsilon \quad \textrm{and} \quad \delta < \lan e_p^*, M e_p \ran \quad \textrm{for every} \; 0 \leq k \leq N.
$$

With these choices, we can easily show that $T$ is a positive contraction on $X$ and that
    $$\lVert (T-M) e_k \rVert < \varepsilon, \quad \lVert (T-M) ^* e_k^* \rVert < \varepsilon \quad \textrm{for every} \; 0 \leq k \leq N.$$

    It remains to show that $T$ does not belong to $\mathcal{F}_{i,j,\eta,p}$. To do so, let $A$ be a positive operator such that $A T = T A$ and $\lan e_p^*, A e_p \ran = 0$. We will prove that $A = 0$. In particular, the condition $\lan e_j^*, A e_i \ran \geq \eta$ won't be fulfilled, and this will yield that $T$ does not belong to $\mathcal{F}_{i,j,\eta,p}$.
\medskip
    
    The equation $A Te_p = T A e_p$ implies that
    \begin{align} \label{eq ep}
    AM e_p + \delta A e_{p+N+1} = MPA e_p + \delta \, \lan e_{p+N+1}^*, A e_p \ran \, u + S A e_p. 
    \end{align}
    We remark that
    $$
    \lan e_{p+N+1}^*, S A e_p \ran = \delta \, \lan e_p^*, A e_p \ran = 0
    $$
    and that $M P A e_p \in E_N$, so we also have
    $$
    \lan e_{p+N+1}^*, MPA e_p \ran = 0.
    $$
    Applying the functional $e_{p+N+1}^*$ to the equation (\ref{eq ep}) one gets that
    \begin{align*}
        \delta \, \lan e_{p+N+1}^* , A e_p \ran &= \lan e_{p+N+1}^*, AM e_p \ran + \delta \, \lan e_{p+N+1}^* , A e_{p+N+1} \ran \\
        &= \displaystyle \sum_{k=0}^N \lan e_k^*, M e_p \ran \, \lan e_{p+N+1}^*, A e_k \ran + \delta \, \lan e_{p+N+1}^*, A e_{p+N+1} \ran
    \end{align*}
    and so
    \begin{align} \label{eqdeltaa}
        (\delta - \lan e_p^*, M e_p \ran) \, \lan e_{p+N+1}^*, A e_p \ran = \displaystyle \sum_{\substack{0 \leq k \leq N  \\ k \ne p}} &\lan e_k^*, M e_p \ran \, \lan e_{p+N+1}^*, A e_k \ran \\ &+ \delta \, \lan e_{p+N+1}^*, A e_{p+N+1} \ran. \notag
    \end{align}
    Using the non-negativity of the coefficients on the right hand side of (\ref{eqdeltaa}), and using the fact that $\delta < \lan e_p^*, M e_p \ran$ and that $\lan e_{p+N+1}^*, A e_p \ran \geq 0$, it follows that 
    \begin{align} \label{eq e_p+n+1}
        \lan e_{p+N+1}^* , A e_{p+N+1} \ran = 0 \quad \textrm{and} \quad \lan e_{p+N+1}^* , A e_k \ran = 0 \quad \textrm{for every} \; 0 \leq k \leq N.
    \end{align}
    Thus one gets
    \begin{align}
    A T e_k &= AM e_k + \delta A e_{k+N+1} \quad \textrm{for every} \; 0 \leq k \leq N,  \\ 
    \textrm{and} \quad  T A e_k &= MPA e_k + S A e_k \quad \textrm{for every} \; 0 \leq k \leq N.
    \end{align}
    The relations $T A e_k = A T e_k$ for $0 \leq k \leq N$ imply that
    \begin{equation}
\left\{
    \begin{array}{lcl}
    MP Ae_0 + S A e_0 & = & AM e_0 + \delta A e_{N+1}\\
    MPA e_1 + S A e_1 & = & AM e_1 + \delta A e_{N+2}\\
    &\vdots&\\
    MP A e_N + S A e_N & = & AM e_N + \delta A e_{2N+1}
    \end{array}
\right.
\end{equation}
and since $P S A e_l = 0$ for every $0 \leq l \leq N$, we obtain that
\begin{equation} \label{systavantlast}
\left\{
    \begin{array}{lcl}
    MP Ae_0  & = & PAM e_0 + \delta PA e_{N+1}\\
    MPA e_1 & = & PAM e_1 + \delta PA e_{N+2}\\
    &\vdots&\\
    MP A e_N  & = & PAM e_N + \delta PA e_{2N+1}.
    \end{array}
\right.
\end{equation}
Finally, rewriting the system (\ref{systavantlast}), we obtain the following system:
\begin{equation} \label{systlast}
\left\{
    \begin{array}{lcl}
    \displaystyle\sum_{k=0}^N \lan e_k^*, PA e_0 \ran M e_k  & = &\displaystyle \sum_{k=0}^N \lan e_k^*, M e_0  \ran P A e_k + \delta PA e_{N+1}\\
    \displaystyle\sum_{k=0}^N \lan e_k^*, PA e_1 \ran M e_k & = & \displaystyle \sum_{k=0}^N \lan e_k^*, M e_1  \ran P A e_k + \delta PA e_{N+2}\\
    &\vdots&\\
    \displaystyle\sum_{k=0}^N \lan e_k^*, PA e_N \ran M e_k  & = & \displaystyle \sum_{k=0}^N \lan e_k^*, M e_N  \ran P A e_k + \delta PA e_{2N+1}.
    \end{array}
\right.
\end{equation}

The system (\ref{systlast}) is equivalent to the following matrix equation

\begin{align} \label{eq systmatrice}
    B C = C B + \delta D,
\end{align}
where
$$
B = \begin{pmatrix} \lan e_0^*, PA e_0 \ran  & \lan e_1^*, PA e_0 \ran & \cdots & \lan e_N^*, PA e_0 \ran \\
\lan e_0^*, PA e_1 \ran  & \lan e_1^*, PA e_1 \ran & \cdots & \lan e_N^*, PA e_1 \ran \\
\vdots  & \vdots & \ddots & \vdots \\
\lan e_0^*, PA e_N \ran & \lan e_1^*, PA e_N \ran & \cdots & \lan e_N^*, PA e_N \ran \\
\end{pmatrix},
$$

$$
C = \begin{pmatrix} \lan e_0^*, M e_0 \ran  & \lan e_1^*, M e_0 \ran & \cdots & \lan e_N^*, M e_0 \ran \\
\lan e_0^*, M e_1 \ran  & \lan e_1^*, M e_1 \ran & \cdots & \lan e_N^*, M e_1 \ran \\
\vdots  & \vdots & \ddots & \vdots \\
\lan e_0^*, M e_N \ran & \lan e_1^*, M e_N \ran & \cdots & \lan e_N^*, M e_N \ran \\
\end{pmatrix},
$$
and
$$
D = \begin{pmatrix} \lan e_0^*, PA e_{N+1} \ran  & \lan e_1^*, PA e_{N+1} \ran & \cdots & \lan e_N^*, PA e_{N+1} \ran \\
\lan e_0^*, PA e_{N+2} \ran  & \lan e_1^*, PA e_{N+2} \ran & \cdots & \lan e_N^*, PA e_{N+2} \ran \\
\vdots  & \vdots & \ddots & \vdots \\
\lan e_0^*, PA e_{2N+1} \ran & \lan e_1^*, PA e_{2N+1} \ran & \cdots & \lan e_N^*, PA e_{2N+1} \ran \\
\end{pmatrix}.
$$
We will now show that $D = 0$ using a similar idea to  \cite[Theorem 2.1]{BDFRZ}.

Since $C$ is a matrix with positive entries, the spectral radius of $C$ and $C^\top$ is a positive eigenvalue of $C$ and $C^\top$ respectively associated to eigenvectors with positive entries. Let $\lambda$ be the spectral radius of $C$ and $C^\top$. Then there exist two vectors $x,y \in \R^{N+1}$ with positive entries such that
\begin{align*}
    C x = \lambda x \quad \textrm{and} \quad y^\top C = \lambda y^\top.
\end{align*}
The equation (\ref{eq systmatrice}) implies that
$$
\lambda Bx = C B x + \delta D x
$$
and 
$$
\lambda y^\top B x = \lambda y^\top B x + \delta y^\top D x,
$$
so that
$$
y^\top D x = 0.
$$
Using the fact that the vectors $x$ and $y$ have positive entries, it follows that $D = 0$.

Now because the vectors $PA e_l$ belong to $E_N$, we have just proved that
\begin{align} \label{eqN+1-2N+1}
    PA e_{k} = 0 \quad \textrm{for every} \; N+1 \leq k \leq 2N+1.
\end{align}

Using the last equation (\ref{eqN+1-2N+1}), the relation $A T e_{p+N+1} = T A e_{p+N+1}$ gives 

$$
\delta A u + \delta A e_{p+2N+2} = S A e_{p+N+1} \,,
$$
so using the fact that
$$
\lan e_k^*, S A e_{p+N+1} \ran = 0 \quad \textrm{for every } \; 0 \leq k \leq N,
$$ 
we get that
\begin{align} \label{eq PAu}
    \delta PAu + \delta PA e_{p+2N+2} =  P S A e_{p+N+1} = 0.
\end{align}
Now recall that $u = e_0 + ... + e_{p+N+1}$. It follows from the non-negativity of the coefficients in the equation (\ref{eq PAu}) that
\begin{align} \label{PAu=0}
    PA e_k = 0 \quad \textrm{for every} \; 0 \leq k \leq 2N+1.
\end{align}

Let
$$
\akl = \lan e_k^*, A e_l \ran \quad \textrm{and} \quad m_{k,l} = \lan e_k^*, M e_l \ran \quad \textrm{for every} \; k,l \geq 0.
$$

From the equations (\ref{eq e_p+n+1}) and (\ref{PAu=0}), we have:
\begin{align}
    &\akl = 0 \quad \textrm{for every} \; 0 \leq k \leq N \; \textrm{and every} \; 0\leq l \leq 2N+1, \\
    &a_{p+N+1,l} = 0 \quad \textrm{for every} \; 0 \leq l \leq N, \\
    &a_{p+N+1,p+N+1} = 0.
\end{align}

An easy computation shows that
$$
\lan e_k^*, ATe_0 \ran = \displaystyle \sum_{l=0}^N a_{k,l} \, m_{l,0} + \delta a_{k,N+1} \quad \textrm{for every} \; k \geq 0 
$$
and that
$$
\lan e_k^*, TA e_0 \ran = \left\{
    \begin{array}{lll}
        \delta \, a_{k-N-1,0} + \delta \, a_{N+p+1,0}  & \mbox{if } N+1 \leq k \leq N+p+1, \\
        \delta \, a_{k-N-1,0} & \mbox{if} \;  N+p+1 < k \leq 2N+p+2,\\
        \delta_{k-2N-p-2} \, a_{k-N-1,0} & \mbox{if} \; k > 2N +p+2.
    \end{array}
\right.
$$
So using the equation 
$$
\lan e_k^*, AT e_0 \ran = \lan e_k^*, TA e_0 \ran \,,
$$
one gets that
$$
\displaystyle \sum_{l=0}^N a_{k,l} \, m_{l,0} + \delta \, a_{k,N+1} = 0 \quad \mbox{for every} \; N+1 \leq k \leq N+p+1,
$$
and using the non-negativity of the coefficients and the fact that the coefficients $m_{l,0}$ are positive, we obtain that
\begin{align}
a_{k,l} = 0 \quad \mbox{for every} \; 0 \leq l \leq N+1 \; \textrm{and every} \, N+1 \leq k \leq N+p+1.
\end{align}
Proceeding by induction, we easily get that
\begin{align}
a_{k,l} = 0 \quad \mbox{for every} \; 0 \leq l \leq N+1 \; \textrm{and every} \; k \geq 0,
\end{align}
so in particular we have
\begin{align} \label{eq Ael}
    A e_l = 0 \quad \mbox{for every} \; 0 \leq l \leq N+1.
\end{align}
Now because we have
\begin{align} \label{eq TAEl0}
    T A e_l = 0 \quad \textrm{for every} \; 0 \leq l \leq N+1
\end{align}
and 
\begin{align} \label{relation ATel}
    A T e_l = A M P e_l + \delta \lan e_{N+p+1}^*, e_l \ran A u + \delta A e_{l+N+1} \quad \textrm{for every} \; 0 \leq l \leq N+1,
\end{align}
the equation (\ref{eq TAEl0}) and the positivity of the vectors in the relation (\ref{relation ATel}) give us
\begin{align} \label{eq e_l+N+1}
    A e_{l+N+1} = 0 \quad \textrm{for every} \; 0 \leq l \leq N+1.
\end{align}

 If we put together the equations (\ref{eq Ael}) and (\ref{eq e_l+N+1}), we obtain that 
$$
A e_l = 0 \quad \mbox{for every} \; 0 \leq l \leq 2N+2.
$$

Proceeding by induction, we easily obtain that $A e_k = 0$ for every $k \geq 0$ and thus $A=0$.

This concludes the proof of Proposition \ref{propdens}.
\end{proof}

The proof of Theorem \ref{thaabreflexivesec4} immediately follows from Fact \ref{factuno} and Proposition \ref{propdens}.
\epf

\smallskip

Finally, we extend Corollary \ref{corqnubcondsec4} to Banach spaces with a basis. Observe that we do not require here that the basis be monotone, and hence Theorem \ref{thaabreflexivesec4} does not apply.

\smallskip

\bpr \label{propqnbasisonly}
Let $X$ be a Banach space with a basis $(e_n)_{n \geq 0}$. A typical $T \in (\ppx, \emph{\sot})$ (resp. $T \in (\ppx, \emph{\sote})$ when $X^*$ is separable) is not quasinilpotent at any non-zero positive vector of $X$.
\epr

\smallskip

\bpf
Consider the following set 
$$\mathcal{A} := \{ T \in \ppx : \exists \, y \in X, \, y \geq 0, \, y \ne 0 \; \textrm{such that} \; T \; \textrm{is quasinilpotent at} \; y  \}.$$
By Lemma \ref{lemqn}, we have
$$
\mathcal{A} \subseteq \displaystyle \bigcup_{j \geq 0} \mathcal{F}_j \,,
$$
with
$$
\mathcal{F}_j := \{ T \in \ppx : \, \lan e_j^*, T e_j \ran = 0 \} \quad \textrm{for every} \; j \geq 0.
$$
Each $\mathcal{F}_j$ is \sot-closed in $\ppx$ and hence \sote-closed in $\ppx$. Let us now prove that each $\mathcal{F}_j$ has empty interior in $(\ppx, \sote)$.

To do so, let $\varepsilon >0$, let $T \in \ppx$ with $\lVert T \rVert < 1$ and let $x_1,...,x_n \in X$ and $y_1^*,...,y_n^* \in X^*$. We have to find a positive contraction $S $ in the set $\ppx \setminus \mathcal{F}_{j}$ such that
\begin{equation} \label{eqqnsec4}
    \displaystyle \max_{1 \leq l \leq n} \max\{ \lVert (T-S) \,x_l \rVert, \lVert (T-S)^* y_l^* \rVert \} < \varepsilon.
\end{equation}
Consider the positive operator $S_\delta$ defined by $S_\delta(x) = T x + \delta \, \lan e_j^*, x \ran \, e_j$ for every $x \in X$, where $\delta$ is a positive number that we will define later on. We have that 
$$
\lan e_j^*, S_\delta \, e_j \ran \geq \delta > 0.
$$
For every $x \in X$, we have that
$$
\lVert S_\delta \, x \rVert \leq (\lVert T \rVert + \delta \lVert e_j^* \rVert \lVert e_j \rVert) \lVert x \rVert
$$
and for every $1 \leq l \leq n$, we have that
$$
\lVert (T - S_\delta) \, x_l \rVert \leq \delta \lVert x_l \rVert \lVert e_j^* \rVert \lVert e_j \rVert \quad \textrm{and} \quad \lVert (T - S_\delta)^* y_l^* \rVert \leq \delta \lVert y_l^* \rVert \lVert e_j^* \rVert \lVert e_j \rVert.
$$
If we choose $\delta > 0$ such that
$$
\delta < \frac{1- \lVert T \rVert}{\lVert e_j^* \rVert \lVert e_j \rVert}, \quad \delta \lVert x_l \rVert \lVert e_j^* \rVert \lVert e_j \rVert < \varepsilon \quad \textrm{and} \quad \delta \lVert y_l^* \rVert \lVert e_j^* \rVert \lVert e_j \rVert < \varepsilon \quad \textrm{for every} \; 1 \leq l \leq n,
$$
then the operator $S_\delta$ is a positive contraction of $\ppx \setminus \mathcal{F}_j$ satisfying (\ref{eqqnsec4}). This proves that each $\ppx \setminus \mathcal{F}_j$ is $\sote$-dense in $\ppx$ and this concludes the proof of Proposition \ref{propqnbasisonly}.
\epf


\section{Further remarks and questions}\label{Questions}
We end this article with some comments and open questions in relation to our previous results.

\par\smallskip
The first natural open question is of course the following. 

\smallskip

\begin{question}\label{Question 1}
  If $X = \ell_q$ with $1 < q \ne 2 < \infty$, does a typical $T \in (\ppx, \sot)$ or $T \in (\ppx,\sote) $ have a non-trivial invariant subspace? 
\end{question}

\smallskip

By \cite[Corollary 5.3]{EM}, the point spectrum of a typical contraction $T \in (\mathcal{B}_1(\ell_2), \sot)$ is the open unit disk $\D$. This comes from the fact that an \sot-typical contraction on $\ell_2$ is unitarily equivalent to the infinite-dimensional backward unilateral shift operator on $\ell_2(\Z_+ \times \Z_+)$. The proof uses first the fact that a typical contraction on $\ell_2$ is a co-isometry for the \sot\, topology. Since this is no longer the case for an \sot-typical positive contraction on $\ell_2$ by Proposition \ref{propcoisolq}, the proof given in \cite{EM} does not work for positive contractions. So the following question is still open. 

\smallskip

\begin{question}\label{Question 2}
 Is it still true that the point spectrum of an \sot-typical positive contraction on $\ell_2$ is $\D$?
\end{question}

\smallskip

The third question is motivated by Proposition \ref{prophypercycl}.

\smallskip

\begin{question}
    If $X = \ell_q$ with $1 < q < 2$, is it true that a typical $T \in (\ppx,\sot)$ is such that $(2T)^*$ is hypercyclic? 
\end{question}

\smallskip

Lemma \ref{lemapprox} requires $X$ to have a monotone basis in order to be able to say that if $T_0$ is a positive operator on $X$ such that $\lVert T_0 \rVert < 1$, then $\lVert P_N T_0 P_N \rVert < 1$ for every $N \geq 0$. This lemma was useful to prove Theorem \ref{thaabreflexivesec4}. We thus have the following open question.

\smallskip

\begin{question}
    Can Lemma \ref{lemapprox} be generalized to Banach spaces admitting a basis which is not necessarily monotone? 
\end{question}

\smallskip

Theorem \ref{thaabreflexivesec4} applies to the \sote\, topology. Indeed, the proof of Fact \ref{factuno} uses the \sote\, topology to prove that each set $\mathcal{F}_{i,j,\eta,p}$ is closed is $\ppx$. Since these sets are not necessarily \sot-closed, the following question is natural.

\smallskip

\begin{question}
    Is the analogue of Theorem \ref{thaabreflexivesec4} still true for the \sot\, topology?
\end{question}

\smallskip

It is proved in \cite[Theorem 7.5]{GMM1} that a typical contraction $T \in (\mathcal{B}_1(\ell_2),\sote)$ does not commute with any non-zero compact operator on $\ell_2$. Since the proof uses unitary equivalence of operators, it does not extend to positive contractions. Thus, the following question is open.

\smallskip

\begin{question}
    Does a typical $T \in (\ppx, \sote)$ commute with a non-zero compact operator if $X = \ell_2$? And what about a typical $T \in (\ppx, \sot)$?
\end{question}

\smallskip

Finally, we have the following question. A positive answer to it would enlighten the situation very much.

\smallskip

\begin{question}
    If $X = \ell_q$ with $1 < q \leq 2$, are the \sot\, and the \sote \,topologies similar on $\ppx $?
\end{question}

\medskip

\textbf{Acknowledgments.} I would like to thank Sophie Grivaux for
helpful discussions and suggestions.


\begin{bibdiv}
  \begin{biblist}
  

\bib{AAB}{article}{
   author={Abramovich, Y. A.},
   author={Aliprantis , C. D.},
   author={Burkinshaw, O.},
   title={Invariant Subspaces for Positive Operators Acting on a Banach Space with Basis},
   journal={Proceedings of the American Mathematical Society},
   volume={123},
   date={1995},
   number={6},
   pages={1773-1777},
}

\bib{AAB2}{article}{
   author={Abramovich, Y. A.},
   author={Aliprantis , C. D.},
   author={Burkinshaw, O.},
   title={Invariant Subspaces of operators on $\ell_p$-spaces},
   journal={Journal of Functional Analysis},
   volume={115},
   date={1993},
   pages={418-424},
}



\bib{AM}{article}{
   author={Ambrozie, C\u{a}lin},
   author={M\"{u}ller, Vladim\'{\i}r},
   title={Invariant subspaces for polynomially bounded operators},
   journal={J. Funct. Anal.},
   volume={213},
   date={2004},
   number={2},
   pages={321--345},
}



\bib{BM}{book}{
   author={Bayart, Fr\'{e}d\'{e}ric},
   author={Matheron, \'{E}tienne},
   title={Dynamics of linear operators},
   series={Cambridge Tracts in Mathematics},
   volume={179},
   publisher={Cambridge University Press, Cambridge},
   date={2009},
   pages={xiv+337},
   }


\bib{BDFRZ}{article}{
    AUTHOR = {Bračič, Janko},
    author={Drnovšek, Roman},
    author={Farforovskaya, Yuliya B.},
    author ={L. Rabkin, Evgueniy},
    author={Zemánek, Jaroslav},
     TITLE = {On positive commutators},
   JOURNAL = {Positivity},
    VOLUME = {14},
      YEAR = {2010},
     PAGES = {431--439},
}


\bib{BCP2}{article}{
    AUTHOR = {Brown, Scott W.},
    author={Chevreau, Bernard},
    author={Pearcy, Carl},
     TITLE = {On the structure of contraction operators. {II}},
   JOURNAL = {J. Funct. Anal.},
    VOLUME = {76},
      YEAR = {1988},
    NUMBER = {1},
     PAGES = {30--55},
}


\bib{CE}{article}{
   author={Chalendar, Isabelle},
   author={Esterle, Jean},
   title={Le probl\`eme du sous-espace invariant},
   conference={
      title={Development of mathematics 1950--2000},
   },
   book={
      publisher={Birkh\"{a}user, Basel},
   },
   date={2000},
   pages={235--267},
}


\bib{CP}{book}{
   author={Chalendar, Isabelle},
   author={Partington, Jonathan R.},
   title={Modern approaches to the invariant-subspace problem},
   series={Cambridge Tracts in Mathematics},
   volume={188},
   publisher={Cambridge University Press, Cambridge},
   date={2011},
   pages={xii+285},
}

\bib{Con}{book}{
    AUTHOR = {Conway, John B.},
     TITLE = {A course in functional analysis},
    SERIES = {Graduate Texts in Mathematics},
    VOLUME = {96},
   EDITION = {second edition},
 PUBLISHER = {Springer-Verlag, New York},
      YEAR = {2007},
     PAGES = {xvi+399},
}


\bib{E}{article}{
   author={Eisner, Tanja},
   title={A ``typical'' contraction is unitary},
   journal={Enseign. Math. (2)},
   volume={56},
   date={2010},
   number={3-4},
   pages={403--410},
}

\bib{EM}{article}{
   author={Eisner, Tanja},
   author={M\'{a}trai, Tam\'{a}s},
   title={On typical properties of Hilbert space operators},
   journal={Israel J. Math.},
   volume={195},
   date={2013},
   number={1},
   pages={247--281},
}


\bib{En}{article}{
   author={Enflo, Per},
   title={On the invariant subspace problem for Banach spaces},
   journal={Acta Math.},
   volume={158},
   date={1987},
   number={3-4},
   pages={213--313},
   
}


\bib{GMM}{article}{
    author={Grivaux, S.},
   author={Matheron, \'{E}.},
   author={Menet, Q.},
   title={Linear dynamical systems on Hilbert spaces: typical properties and
   explicit examples},
   journal={Mem. Amer. Math. Soc.},
   volume={269},
   date={2021},
   number={1315},
   pages={v+147},
}

\bib{GMM1}{article}{
   author={Grivaux, Sophie},
   author={Matheron, \'{E}tienne},
   author={Menet, Quentin},
   title={Does a typical $\ell_p$-space contraction have a non-trivial
   invariant subspace?},
   journal={Trans. Amer. Math. Soc.},
   volume={374},
   date={2021},
   number={10},
   pages={7359--7410},

}

\bib{GMM2}{article}{
   author={Grivaux, Sophie},
   author={Matheron, \'{E}tienne},
   author={Menet, Quentin},
   title={Generic properties of $\ell_p$-contractions and similar operator topologies},
   journal={preprint available at https://arxiv.org/abs/2207.07938},
   date={2022},   

}
		





\bib{GEP}{book}{
   author={Grosse-Erdmann, Karl-G.},
   author={Peris Manguillot, Alfredo},
   title={Linear chaos},
   series={Universitext},
  publisher={Springer, London},
   date={2011},
   
}


\bib{Ke}{book}{
   author={Kechris, Alexander S.},
   title={Classical descriptive set theory},
   series={Graduate Texts in Mathematics},
   volume={156},
   publisher={Springer-Verlag, New York},
   date={1995},
   pages={xviii+402},
   isbn={0-387-94374-9},
   review={\MR{1321597}},
   doi={10.1007/978-1-4612-4190-4},
}



\bib{LT}{book}{
   author={Lindenstrauss, Joram},
   author={Tzafriri, Lior},
   title={Classical Banach spaces. I},
   note={Sequence spaces;
   Ergebnisse der Mathematik und ihrer Grenzgebiete, Vol. 92},
   publisher={Springer-Verlag, Berlin-New York},
   date={1977},
   pages={xiii+188},
}


\bib{L}{article}{
   author={Lomonosov, V. I.},
   title={Invariant subspaces of the family of operators that commute with a
   completely continuous operator},
   language={Russian},
   journal={Funkcional. Anal. i Prilo\v{z}en.},
   volume={7},
   date={1973},
   number={3},
   pages={55--56},
}



\bib{Mul}{book}{
   author={Müller, Vladimir},
   title={Spectral Theory of Linear Operators},
   series={Operator Theory: Advances and Applications},
   volume={139},
   edition={2},
   publisher={Birkhäuser Basel},
   date={2007},
   pages={xvii+439},
}

		




\bib{RR}{book}{
   author={Radjavi, Heydar},
   author={Rosenthal, Peter},
   title={Invariant subspaces},
   edition={2},
   publisher={Dover Publications, Inc., Mineola, NY},
   date={2003},
   pages={xii+248},
}

\bib{RT}{article}{
   author={Radjavi, Heydar},
   author = {Troitsky, Vladimir},
   title={Invariant sublattices},
   journal={Ilinois J. Math.},
   volume={52},
   date={2008},
   number={2},
   pages={437--462},
}



\bib{R2}{article}{
   author={Read, C. J.},
   title={A solution to the invariant subspace problem on the space $l_1$},
   journal={Bull. London Math. Soc.},
   volume={17},
   date={1985},
   number={4},
   pages={305--317},
   
}
		


\bib{R3}{article}{
   author={Read, C. J.},
   title={The invariant subspace problem on some Banach spaces with
   separable dual},
   journal={Proc. London Math. Soc. (3)},
   volume={58},
   date={1989},
   number={3},
   pages={583--607},
   
}		

		

\bib{Troi}{article}{
   author={Troitsky, V.},
   title={On the modulus of C. J. Read's Operator},
   journal={Positivity},
   volume={2},
   date={1998},
   pages={257--264},
}



  \end{biblist}
\end{bibdiv}
\end{document}